\documentclass[10pt,letter]{article}

\usepackage{lipsum}

\usepackage{tikz-cd}
\usepackage{amsmath}
\usepackage{amssymb}
\usepackage{amsthm, thmtools}
\usepackage{mathtools}
\usepackage{mathrsfs}

\usepackage{bbm}
\usepackage{bm}
\usepackage{MnSymbol}   
\usepackage{cases}
\usepackage[normalem]{ulem}

\usepackage[width=.80\textwidth,font=footnotesize]{caption}
\usepackage[letterpaper, margin=1.25in]{geometry}
\usepackage[title]{appendix}
\usepackage{fancyhdr}
\usepackage{caption}
\usepackage{subcaption}

\usepackage[shortlabels]{enumitem}
\usepackage[ddmmyyyy]{datetime}
\usepackage{hhline}

\usepackage{listings}
\newcommand{\N}{\mathbb{N}}

\newcommand{\R}{\mathbb{R}}

\newcommand{\pr}{\mathbb{P}}

\DeclarePairedDelimiter\abs{|}{|}
\DeclarePairedDelimiterX{\inner}[2]{\langle}{\rangle}{#1,#2}
\DeclarePairedDelimiterX{\norm}[1]{\|}{\|}{#1}

\DeclareMathOperator{\lip}{Lip}
\DeclareMathOperator{\spec}{Spec}
\DeclareMathOperator{\Tr}{Tr}

\DeclareMathOperator{\supp}{supp}

\DeclareMathOperator{\ExpOp}{\mathbb{E}}

\newtheorem{innercustomthm}{Theorem}
\newenvironment{customthm}[1]
  {\renewcommand\theinnercustomthm{#1}\innercustomthm}
  {\endinnercustomthm}

\newtheorem{theorem}{Theorem}[section]
\newtheorem*{theorem*}{Theorem}
\newtheorem{lemma}[theorem]{Lemma}	\newtheorem*{lemma*}{Lemma}
\newtheorem{proposition}[theorem]{Proposition}
\newtheorem*{cor*}{Corollary}

\theoremstyle{definition}	

\newtheorem{assumption}{Assumption}

\theoremstyle{definition}

\theoremstyle{plain}		
\theoremstyle{definition}	\newtheorem{remark}[theorem]{Remark}

\theoremstyle{definition}	\newtheorem{example}{Example}[section]

\theoremstyle{definition}

\begin{document}
\begin{center}
{\Large Quasi-Ergodicity of Transient Patterns in Stochastic Reaction-Diffusion Equations}
\end{center}

\begin{center}
\begin{minipage}[t]{.75\textwidth}
\begin{center}
Zachary P.~Adams \\
\footnotesize{Freie Universit{\"a}t Berlin, Scads.AI Leipzig,\\ \&~Max Planck Institute for Mathematics in the Sciences}\\
zachary.adams@mis.mpg.de\\
$\,$\\
\hfill
\end{center}
\end{minipage}
\end{center}
\begin{center}
\today 
\end{center}
\vspace{.5cm}
\setlength{\unitlength}{1in}

\vskip.15in

{
\centerline
{\large \bfseries \scshape Abstract}
We study transient patterns appearing in a class of SPDE using the framework of quasi-stationary and quasi-ergodic measures. 
In particular, we prove the existence and uniqueness of quasi-stationary and quasi-ergodic measures for a class of reaction-diffusion systems perturbed by additive cylindrical noise. 
We obtain convergence results in $L^2$ and almost surely, and demonstrate an exponential rate of convergence to the quasi-stationary measure in an $L^2$ norm. 
These results allow us to qualitatively characterize the behaviour of these systems in neighbourhoods of an invariant manifold of the corresponding deterministic systems at some large time $t>0$, conditioned on remaining in the neighbourhood up to time $t$. 
The approach we take here is based on spectral gap conditions, and is not restricted to the small noise regime.
}\\


\noindent{\bfseries \emph{Published version}}~\emph{available at https://doi.org/10.1214/24-EJP1130}
\hfill

\section{Introduction} 
The purpose of this article is to provide a framework with which to discuss the stability and long time behaviour of stochastic perturbations of spatiotemporal patterns appearing in reaction-diffusion equations. 
This topic has recently been of great interest to physicists, biologists, and mathematicians. 
For instance, travelling waves and stochastic perturbations thereof are important features of models from neuroscience \cite{C05,E72i,IM16,MR81}, population genetics \cite{F58,MMQ11,MMR21}, ecology \cite{RMF08,SS08}, and many other disciplines. 
Meanwhile, understanding stochastically perturbed spiral waves in excitable media has lead to insights on cardiac arrhythmias and how to treat them \cite{BE93,F02,L21}. 

Stochastic perturbations of a pattern in a reaction-diffusion system usually destroy the pattern at some finite time. 
Hence, while the unperturbed pattern may be stable, the perturbed pattern is often only \emph{metastable}. 
In recent decades, there has been much interest in characterizing metastable behaviour in a rigorous mathematical framework, for instance in \cite{BW12,EGK20,IM16,KS17,L16,LS16,S13}. 
In this paper, we approach the study of metastable patterns in stochastic reaction-diffusion systems using the theory of \emph{quasi-stationary}~and \emph{quasi-ergodic measures}, as defined for instance in \cite{CV16}. 
As described more precisely below, these measures characterize the long time behaviour of a metastable pattern prior to its destruction. 
The connection between metastability and quasi-ergodicity has been studied previously, for instance in discrete and one-dimensional settings in \cite{BGM20,BG16,GMV20}~and \cite{H22,JQSY22}, respectively. 
To the author's knowledge, the results of this paper are the first contribution studying metastability via quasi-stationary measures in an infinite dimensional setting. 

The remainder of this document is structured as follows. 
In Section \ref{sec:Overview}, we outline sufficient conditions for the existence and uniqueness of quasi-stationary and quasi-ergodic measures in SPDE, while Section \ref{sec:Literature}~presents an incomplete review of relevant literature. 
In Section \ref{sec:conditionedSemigroups}, we provide results on quasi-ergodicity of general continuous Markov processes in separable Hilbert space. 
As we do not assume a modified Doeblin condition as in \cite{CV16}, nor do we assume the existence of bounded integral kernels for the Markov semigroups which we study as in \cite{CLMR21,HK18,P85,P90,ZLS14}, we obtain convergence results that hold in an $L^2$ and almost sure sense, rather than uniformly. 
In Section \ref{sec:Examples}, we prove that a large class of semilinear SPDEs satisfy the hypotheses of Section \ref{sec:conditionedSemigroups}, and therefore admit unique quasi-stationary and quasi-ergodic measures. 
In Section \ref{sec:Applications}, we conclude the paper with a cursory discussion of how quasi-ergodic measures relate to metastable patterns in SPDEs.

\subsection{Setup and Results}
\label{sec:Overview}
Let $O\subset\R^d$ be a spatial domain and let $H$ be a separable Hilbert space of functions $f:O\rightarrow\R^n$ for some $n\in\N$. 
Consider the following evolution equation on $H$,
\begin{equation}
\label{eq:PDE} 
\partial_t x\,=\,Lx + N(x), 
\end{equation}
where $L$ and $N$ are linear and nonlinear operators, respectively, on $H$. 
We impose the following conditions on \eqref{eq:PDE}, which we verify for the case where $L=\Delta$ is the Laplace operator and $N$ is a polynomial in Example \ref{EX}~below.  
Under these assumptions, unique mild solutions to \eqref{eq:PDE}~are known to exist, for instance using fixed point arguments similar to those discussed in Section \ref{sec:Examples}~below. 

\begin{assumption}
\label{assn:PDE}
The PDE \eqref{eq:PDE}~satisfies the following. 
\begin{enumerate}[(a)]
\item $N$ is a nonlinearity defined on a set $D(N)$ that is densely and continuously embedded in $H$, $D(N)$ possesses its own Banach space structure, and $N:D(N)\rightarrow H$ is locally Lipschitz and twice continuously Fr{\'e}chet differentiable as an operator on $D(N)$. 
\item $L:D(L)\subset H\rightarrow H$ generates a strongly continuous semigroup $(\Lambda_t)_{t\ge0}$ on $H$ such that $\Lambda_t(H)\subset D(N)$ for all $t>0$, and $(\Lambda_t)_{t\ge0}$ restricted to $D(N)$ is a strongly continuous semigroup. 
\item There exists $\omega>0$ such that $\norm*{\Lambda_t}_H\le e^{-\omega t}$ and $\norm*{\Lambda_t}_{D(N)}\le\,e^{-\omega t}$. 
\item $\Lambda_t:H\rightarrow D(N)$ is compact, and $\norm*{\Lambda_t}_{H\rightarrow D(N)}$ is locally square integrable in $t>0$.  
\item There is a stable normally hyperbolic invariant manifold $\Gamma$ of \eqref{eq:PDE}, in the sense of \cite[Condition (H3)]{BLZ98}. 
\end{enumerate}
\end{assumption}

Our primary interest is in how the solutions to \eqref{eq:PDE}~near the invariant manifold $\Gamma$ behave under perturbations by noise. 
Specifically, consider the following stochastic perturbation of \eqref{eq:PDE}, 
\begin{equation}
\label{eq:SPDE}
dX\,=\,(LX + N(X))\,dt + \sigma B\,dW. 
\end{equation}
The noise amplitude is $\sigma\ge0$, $B^2$ is a symmetric positive bounded linear operator with bounded inverse, and $W=(W_t)_{t\ge0}$ is a cylindrical Wiener process on $H$. 

For each $x\in H$, all stochastic objects are considered over a fixed ambient probability space $(\Omega,\mathcal{F},\pr_x)$, where $\pr_x$ is the probability measure associated with \eqref{eq:SPDE}~with initial condition $X_0=x$. 
Let $\ExpOp_x$ denote expectation with respect to $\pr_x$. 
Let $(\mathcal{F}_t)_{t\ge0}$ be the filtration associated with $(W_t)_{t\ge0}$. 
For any probability measure $\nu$ on $H$ let 
\[
\pr_\nu[\cdot]\,\coloneqq\,\int_H\pr_x[\cdot]\,\nu(dx), 
\]
and let $\ExpOp_\nu$ be expectation with respect to $\pr_\nu$. 
For any measure $\nu$ on $H$ and $\nu$-measurable function $f:H\rightarrow\R$, we denote $\nu(f)\coloneqq\int f(x)\,\nu(dx)$. 

We now provide sufficient conditions for the existence and uniqueness of mild solutions to \eqref{eq:SPDE}. 
By a \emph{mild solution}~to \eqref{eq:SPDE}, we mean a continuous Markov process $(X_t)_{t\ge0}$ satisfying 
\[
X_t\,=\,\Lambda_tX_0 + \int_0^t \Lambda_{t-s}N(X_s)\,ds + \sigma\int_0^t\Lambda_{t-s}B\,dW_s,\qquad X_0\in H,\,\,t>0. 
\]

\begin{assumption}
\label{assn:SPDE}
The SPDE \eqref{eq:SPDE}~satisfies the following. 
\begin{enumerate}[(a)] 
\item $Y_t\coloneqq\int_0^t\Lambda_{t-s}B\,dW_s$ is a well-defined continuous $D(N)$-valued process. 
\item For all $t>0$, 
\begin{equation}
\label{eq:Tr}
\Tr\left(\int_0^t \Lambda_s B^2\Lambda_s^*\,ds\right)\,<\,\infty. 
\end{equation}
\end{enumerate}
\end{assumption}

It should be emphasized that while we allow for non-trace class noise, our assumptions guarantee that the SPDE \eqref{eq:SPDE}~is non-singular, so that we need not be concerned with the renormalization theory developed \emph{e.g.}~in \cite{DPD03,H13,H14}. 

In Theorem \ref{thm:Solution}~below, we see that Assumptions \ref{assn:PDE}~\&~\ref{assn:SPDE}~imply the existence and uniqueness of local in time $D(N)$-valued mild solutions to \eqref{eq:SPDE}. 
We let $\Gamma_\delta$ be a $\delta$-neighbourhood of $\Gamma$ defined in the $D(N)$ topology, and denote by $N_\delta$ a $C^2$ cutoff of $N$ supported on a small neighbourhood of $\Gamma_\delta$. 
The construction of $N_\delta$ is specified more carefully in Section \ref{sec:Examples}. 
In Theorem \ref{thm:CutoffInvariant}, we see Assumptions \ref{assn:PDE}~\&~\ref{assn:SPDE}~also imply that \eqref{eq:SPDE}~with $N$ cutoff to $\Gamma_\delta$ possesses a unique invariant measure, henceforth denoted $\mu$, such that $\mu(D(N))=1$. 

Following these considerations, we make the following additional assumption.  

\begin{assumption}
\label{assn:Supp}
Assumptions \ref{assn:PDE}~\&~\ref{assn:SPDE}~hold, so that there exists a unique invariant measure $\mu$ of \eqref{eq:SPDE}~with $N$ replaced by a smooth cutoff $N_\delta$ supported on a neighbourhood of $\Gamma_\delta$, as proven in Theorem \ref{thm:CutoffInvariant}~below. 
We morevoer assume that there exists $\delta>0$ such that the following hold. 
\begin{enumerate}[(a)]
\item $\Gamma_\delta$ is contained in the support of $\mu$, so that $\mu(\Gamma_\delta)>0$. 
\item Denoting the Lipschitz constant of $N_\delta$ by $\kappa$, we have $\kappa\in(0,\omega)$. 
\end{enumerate}
\end{assumption}

\begin{example}
\label{EX}
The principal example we have in mind is as follows. 
Let $O\subset\R^d$ be a sufficiently regular bounded spatial domain and set $H\coloneqq L^2(O,\R^n)$. 
Let $L=d\Delta$ with Dirichlet boundaries for some $d>0$, or $L=d\Delta-a$ for constants $a,d>0$ with periodic boundaries. 
For some vector of polynomials $p:\R^n\rightarrow\R^n$, define 
\[
N(f)(\xi)\,=\,p(f(\xi))\quad\text{ for } \,\,f\in D(N)\,\coloneqq\,C_b(O,\R^n),\,\,\,\xi\in O. 
\]
It is straightforward to check that this example satisfies Assumption \ref{assn:PDE}(a), (b), (c), (d), due to the compactness of $\Lambda_t:H\rightarrow D(N)$ and the $L^2\rightarrow L^\infty$ ultracontractivity of $(\Lambda_t)_{t\ge0}$ \cite[Theorem 9.3]{D84}. 
In this example, we suppose that \eqref{eq:PDE}~admits a stable travelling wave solution, as in the FitzHugh-Nagumo equation on $O=[0,L]$ with periodic boundaries \cite{AK15,FH61,KSS97}. 
Specifically, we suppose that there exists $\hat{x}\in C_b(O;\R^n)$ and $c\in\R^d$ such that 
\[
x_t(\xi)\,=\,\hat{x}(\xi-ct),\qquad \xi\,\in\,O,\,\,\,t\,>\,0, 
\]
is a solution of \eqref{eq:PDE}, and $\Gamma\coloneqq\left\{\hat{x}(\cdot-ct)\,:\,t\in\R\right\}$ is a nonlinearly stable invariant manifold of \eqref{eq:PDE}. 
Hence there exists $\delta>0$ such that if 
\[
\sup_{\xi\in O}\norm*{x_0(\xi)-\hat{x}(\xi-cs)}_{\R^d}\,\le\,\delta 
\]
for some $s\in\R$ and $(x_t)_{t\ge0}$ is the solution to \eqref{eq:PDE}~with initial condition $x_0$, then 
\[
\lim_{t\rightarrow\infty}d(x_t,\Gamma)\,=\,0, 
\]
where $d(x,A)$ is the distance between a point $x$ and the closest point in a set $A$. 
So long as we work in one spatial dimension, we may take $B=I$ and verify Assumption \ref{assn:SPDE}~as in \cite[Theorem 11.3.1]{DPZ96}. 
If we take $B$ such that $\Tr B<\infty$, then we may take $O\subset\R^d$ for $d>1$, and verify that Assumption \ref{assn:SPDE}~is still satisfied. 
In this example, by taking sufficiently small $\delta$ we can verify Assumption \ref{assn:Supp}~so long as either of the constants $a,d$ in $L=d\Delta-a$ are large enough for Assumption \ref{assn:Supp}(b) to hold. 
\end{example}

Now, we return to the abstract SPDE \eqref{eq:SPDE}. 
Fix $\delta>0$ such that Assumption \ref{assn:Supp}~holds. 
As we are interested in the dynamics of \eqref{eq:SPDE}~in $\Gamma_\delta$, it is natural to define 
\[
\tau\,\coloneqq\,\inf\left\{t>0\,:\,X_t\in H/\Gamma_\delta\right\}. 
\]
We say that a solution of \eqref{eq:SPDE}~is \emph{$\Gamma$-like}~(with tolerance $\delta>0$) at time $t>0$ if $X_t\in\Gamma_\delta$. 
We would like to be able to answer the following questions. 
\begin{enumerate}
\item[{\bf Q1.}] (Stability) If $X_0\in\Gamma$, how long does a $\Gamma$-like solution persist? 
\item[{\bf Q2.}] (Similarity) If $X_t$ is a $\Gamma$-like solution over a (random) time interval $[0,\tau)$, how does its behaviour differ from the behaviour of solutions to \eqref{eq:PDE}~in $\Gamma$ over this time interval? 
\item[{\bf Q3.}] (Ergodicity) When, and how, can we understand the ``average'' behaviour of a $\Gamma$-like solution? 
\end{enumerate}

A brief, non-exhaustive survey of some of the literature related to these questions is presented in Section \ref{sec:Literature}. 
In this paper, we address these questions -- primarily Q3 -- using the theory of quasi-stationary and quasi-ergodic measures. 
We obtain the following result on the existence and uniqueness of quasi-stationary and quasi-ergodic measures for \eqref{eq:SPDE}. 

\begin{customthm}{A}
Let Assumption \ref{assn:Supp}~hold. 
Then, there exist unique positive $\varphi,\varphi^*\in L^2(\Gamma_\delta,\mu)$, given by the principal eigenfunctions of the sub-Markov semigroup of \eqref{eq:SPDE}~with killing at $H/\Gamma_\delta$, such that the following statements hold. 
\begin{enumerate}[(i)]
\item $\alpha(dx)\coloneqq\varphi^*(x)\mu(dx)$ is the unique probability measure on $\Gamma_\delta$ satisfying  
\begin{equation}
\label{eq:QSD0}
\pr_\alpha\left[X_t\in A\,|\,t<\tau\right]\,=\,\alpha(A) 
\end{equation}
for all measurable $A\subset\Gamma_\delta$. 
Moreover, for $f\in L^2(\Gamma_\delta,\mu)$, the following convergence result holds in $L^2_\mu$ and $\mu$-almost surely: 
\begin{equation}
\label{eq:QSD0_lim}
\lim_{t\rightarrow\infty}\ExpOp_x\left[f(X_t)\,|\,t<\tau\right]\,=\,\alpha(f). 
\end{equation}
\item $\beta(dx)\coloneqq\varphi(x)\varphi^*(x)\mu(dx)$ is the unique probability measure on $\Gamma_\delta$ such that for any $\epsilon>0$, and bounded measurable $f:H\rightarrow\R$,
\begin{equation}
\label{eq:QED0}
\pr_x\left[\norm*{\frac{1}{t}\int_0^t f(X_s)\,ds - \beta(f)}>\epsilon\,|\,t<\tau\right]\,\xrightarrow[t\rightarrow\infty]{}\,0, 
\end{equation}
the above convergence holding in $L^2_\mu$ and $\mu$-almost surely over $x\in\Gamma_\delta$. 
\end{enumerate}
\end{customthm}

The proof of Theorem A follows from Theorems \ref{thm:quasiErgodic}~\&~\ref{thm:XtPt}~below. 
A probability measure satisfying \eqref{eq:QSD0}~is referred to as a \emph{quasi-stationary measure}, while a probability measure satisfying \eqref{eq:QED0}~is referred to as a \emph{quasi-ergodic measure}. 
As will be seen in Section \ref{sec:conditionedSemigroups}, the functions $\varphi,\varphi^*$ are given by the principal eigenfunctions of the sub-Markov generator of \eqref{eq:SPDE}~in $L^2(\Gamma_\delta,\mu)$. 
We obtain the following result on the rate of convergence to the quasi-stationary measure $\alpha$, the proof of which follows from Theorems \ref{thm:expQSD}~\&~\ref{thm:XtPt}~below. 

\begin{customthm}{B}
Let Assumption \ref{assn:Supp}~hold, and let $\varphi,\varphi^*,\alpha$ be as in Theorem A. 
Then, defining $\alpha^*(dx)\coloneqq\varphi(x)\mu(dx)$, for any $f\in L^2(\Gamma_\delta,\mu)$ there exists $T>0$ and a constant $C_f>0$ such that  
\[
\norm*{\ExpOp_x\left[f(Z_t)\,|\,t<\tau\right]-\alpha(f)}_{L^1_{\alpha^*}}\,\le\,C_fe^{-\gamma t}\qquad\forall t\ge T. 
\]
\end{customthm}

Building on the above results, we prove the existence of a \emph{$Q$-process}~for \eqref{eq:SPDE}~in $\Gamma_\delta$. 
Roughly speaking, this $Q$-process is a homogeneous Markov process in $\Gamma_\delta$ given by $(X_t)_{t\ge0}$ conditioned on never exitting $\Gamma_\delta$. 
In practice, $Q$-processes may be analyzed via Doob's $h$-transform \cite[pp.~296]{RW00}, for instance as seen in \cite[Theorem 6.16(iv)]{CMSM13}, though this is not a technique used in this document. 
Since we do not prove uniform convergence to the quasi-stationary measure, and instead rely on the spectral properties of the sub-Markov generator of \eqref{eq:SPDE}~in $\Gamma_\delta$, we cannot use the theory developed \emph{e.g.}~in \cite{CV16}. 
The proof of the following theorem is provided in Theorems \ref{thm:Qproc}~\&~\ref{thm:XtPt}. 

\begin{customthm}{C}
Let Assumption \ref{assn:Supp}~hold for some fixed $\delta>0$. 
Then, for each $s\ge0$ and $A\in\mathcal{F}_s$, the limit 
\begin{equation}
\label{eq:Qx0}
\mathbb{Q}_x[A]\,\coloneqq\,\lim_{t\rightarrow\infty}\pr_x\left[A\,|\,t<\tau\right],\qquad A\in\mathcal{F}_s,\,\,s\ge0  
\end{equation}
is defined for $\mu$-almost all $x\in\Gamma_\delta$. 
With respect to $(\mathbb{Q}_x)_{x\in \Gamma_\delta}$, the solution process $(X_t)_{t\ge0}$ is a homogeneous Markov process in $\Gamma_\delta$ with unique ergodic measure $\beta$. 
\end{customthm}

Before proving Theorems A, B, \&~C, we discuss the existence of quasi-stationary and quasi-ergodic measures in a more general Hilbert space setting in Section \ref{sec:conditionedSemigroups}, below, and return to the question of quasi-ergodicity of $\Gamma$-like solutions of \eqref{eq:SPDE}~in Section \ref{sec:Examples}.  
First, we review the literature related to the problem at hand.

\subsection{Literature}
\label{sec:Literature} 
As described in Example \ref{EX}, are principally interested in the metastability of travelling waves in stochastic reaction-diffusion equations, a
Previous studies on this topic include Bresslof \&~Weber \cite{BW12}, Eichinger, Gnann, \&~Kuehn \cite{EGK20}, Hamster \&~Hupkes \cite{H20,HH20}, Inglis \&~MacLaurin \cite{IM16}, Kr{\"u}ger \&~Stannat \cite{KS14,KS17}, Lang \cite{L16}, Lang \&~Stannat \cite{LS16}, MacLaurin \cite{M20}, and Stannat \cite{S13}. 
We also expect the theory developed in this paper to apply to spiral waves on bounded discs in $\R^2$ with Dirichlet boundary conditions, such as those studied by Xin \cite{X93}, and with significant modifications to spiral waves on $\R^2$, as studied in PDE by Beyn \&~Lorenz \cite{BL08}, Sandstede, Scheel, \&~Wulff \cite{SSW97}, and in SPDEs by Kuehn, MacLaurin, \&~Zucal \cite{KMZ22}. 

Other studies have addressed metastability of patterns in SPDE using techniques related to the theory of large deviations \cite{FW98}. 
The recent work of Salins, Budhiraja, \&~Dupuis \cite{SBD19}~provides a general overview of large deviation theory for SPDEs, while Barret \cite{B15}~and Berglund \&~Gentz \cite{BG13}~have obtained Eyring-Kramers formula for SPDEs using finite dimensional approximation methods. 
The results of \cite{B15,BG13}~provide quantitative estimates on the exit time up to which a metastable pattern in an SPDE persists. 
However, these exit time estimates are only valid in a small noise regime, making them irrelevant for many applications to the life sciences, where small noise approximations are often bad. 
Many biological systems are composed of populations of interacting particles/agents, and in practice the population size of such systems is rarely large enough for a deterministic limit to be a good approximation of the dynamics. 
The resulting \emph{demographic noise}~-- \emph{i.e.}~noise arising from the fact that the system is composed of finitely many agents -- can have significant effects on the dynamics of a systems. 
See \cite{BG09,CRM16,GV21,RMF08,R12}~for examples from various fields of biology. 


Metastability has been studied using quasi-stationary and quasi-ergodic measures in finite dimensional systems in \cite{BGM20,BG16,H22,JQSY22}. 
In \cite{GMV20}, the hypotheses of \cite{BGM20,BG16}~are checked for a particular example of the Glauber dynamics associated with a two dimensional Ising model with boundary conditions. 
An introduction to quasi-ergodicity can be found in the textbook of Collet \emph{et al.}~\cite{CMSM13}. 
The questions of existence and uniqueness of quasi-stationary and quasi-ergodic measures are discussed in their greatest generality in the work of Champagnat \&~Villemonais \cite{CV16,CV17a,CV17b,CV21}. 
In \cite{CV16}, exponential convergence to quasi-stationary measures in total variation norm is shown to be equivalent to a modified Doeblin condition, while \cite{CV21}~studies Lyapunov conditions for the existence of quasi-stationary and quasi-ergodic measures. 
Castro \emph{et al.}~\cite{CLMR21}, Hening \&~Kolb \cite{HK18}, Hening \emph{et al.}~\cite{H22}, and Ji \emph{et al.}~\cite{JQSY22}, Leli{\`e}vre \emph{et al.}~\cite{LRR21}, Pinsky \cite{P85,P90}, Zhang, Li, \&~Song \cite{ZLS14}, and others have studied quasi-ergodicity by exploiting the spectral properties of the sub-Markov semigroup of a killed Markov process. 
Their arguments are similar in spirit to those of the present document. 
However, \cite{CLMR21,HK18,P85,P90,ZLS14}~work in finite dimensional settings, and implicitly or explicitly assume the existence of a bounded density of the sub-Markov transition kernel with respect to some referene measure, which we do not assume here. 
Meanwhile \cite{H22,JQSY22,LRR21}~make explicit use of the finite dimensional nature and gradient structure of the systems they study. 

To the author's knowledge, there is almost no work on quasi-stationary and quasi-ergodic measures of SPDEs. 
Liu \emph{et al.}~\cite{LRSS21}~prove quasi-stationarity of subcritical superprocesses, the distribution of which is governed by an SPDE \cite{W68}. 
However, the arguments of \cite{LRSS21}~are only at the level of the particle process, and hence do not generalize to other SPDEs which are not dual to a particle system. 
To the author's knowledge, this paper therefore represents the first results on the existence and uniqueness of quasi-stationary and quasi-ergodic measures for general semilinear SPDEs.

\section{General Sub-Markov Semigroups}
\label{sec:conditionedSemigroups}
In this section, we establish conditions for the existence and uniqueness of quasi-stationary and quasi-ergodic measures for a general irreducible continuous Markov process $(Z_t)_{t\ge0}$ on a separable Hilbert space $H$. 
Our strategy is to demonstrate that if the Markov semigroup of $(Z_t)_{t\ge0}$ is irreducible and compact in some topology, then these properties are inherited by the sub-Markov semigroup of $(Z_t)_{t\ge0}$ with killing at the boundary of any bounded connected subdomain of $H$. 
This allows us to conclude that the top of the spectrum of this sub-Markov generator is a simple eigenvalue. 
These considerations are proven in Section \ref{sec:Compact}. 
Then, in Section \ref{sec:QuasiQ}~we prove quasi-ergodicity of the process with killing. 
We emphasize that we make no assumptions on the existence or boundedness of an integral kernel density for our sub-Markov semigroup, nor do we assume any Doeblin type or Lyapunov condition on our Markov process. 
As a consequence we obtain convergence results that hold in an $L^p$ and almost sure sense, rather than uniformly. 

Before proceeding, we introduce the general setup and assumptions of this section. 
In Section \ref{sec:Examples}, we verify that these assumptions hold for SPDEs satisfying Assumption \ref{assn:Supp}. 
Let $(Z_t)_{t\ge0}$ be a continuous $(\mathcal{F}_t)_{t\ge0}$-adapted Markov process on a separable Hilbert space $H$. 
Let the Markov semigroup of $(Z_t)_{t\ge0}$ be defined as 
\[
P_tf(x)\,\coloneqq\,\ExpOp_x\left[f(Z_t)\right],\quad f\in BM(H),\,\,\,x\in H, 
\]
where $BM(H)$ is the set of bounded measurable functions $f:H\rightarrow\R$. 
In our general setting, $(P_t)_{t\ge0}$ is not necessarily strongly continuous on $BM(H)$. 
However, in many cases, one can in fact find an extension of $BM(H)$ onto which $(P_t)_{t\ge0}$ extends to a strongly continuous semigroup. 
This leads to the following assumptions. 

\begin{assumption}
\label{assn:Compact}
There exists an invariant measure $\mu$ of $(P_t)_{t\ge0}$ such that $P_t$ extends to a compact operator on $L^2(H,\mu)$ for each $t>0$. 
\end{assumption}

\begin{assumption}
\label{assn:Irreducible}
The semigroup $(P_t)_{t\ge0}$ is irreducible, in the sense that 
\begin{equation}
\label{eq:Irreducible}
P_t1_F(x)\,>\,0\qquad\text{for all open measurable $F\subset H$, $x\in H$, $t>0$}, 
\end{equation}
and strong Feller, in the sense that it maps bounded measurable functions to bounded continuous functions. 
\end{assumption}

Now, fix a bounded, connected, nonempty, open subset $E$ of $H$ contained in the support of $\mu$. 
For $E\subset \supp\mu$ and $p\ge1$, define the $L^p_\mu$-norm of $f\in BM(E)$ as 
\begin{equation}
\norm*{f}^p_{L^p_\mu}\,\coloneqq\,\int_E \norm*{f(x)}^p\,\mu(dx).  
\end{equation}
We then define $L^p(E,\mu)$ as the closure of $BM(E)$ under the $L^p_\mu$-norm. 
Since we have specified that $E$ is nonempty and open, and that $E\subset\supp\mu$, these spaces are nontrivial. 
Supposing that the initial distribution of $(Z_t)_{t\ge0}$ has support in $E$, define the stopping time 
\[
\tau\,\coloneqq\,\inf\left\{t>0\,:\,Z_t\in H/E\right\}. 
\]
For $t\ge0$ the sub-Markov semigroup of $(Z_t)_{t\ge0}$ killed on $H/E$ is 
\[
Q_tf(x)\,\coloneqq\,\ExpOp_x\left[f(Z_t)1_{\{t<\tau\}}\right],\qquad f\in C_b(E). 
\]

\begin{assumption}
\label{assn:tau}
For $B\subset E$, let $\tau_B\coloneqq\inf\left\{t>0\,:\,Z_t\in H/B\right\}$, so $\tau=\tau_E$. 
We assume that for all nonempty open $B\subset E$, $x\in E$ and $t>0$, we have $\pr_x\left[\tau_B<\infty\right]=1$ and $\pr_x\left[t<\tau_B\right]>0$. 
\end{assumption}

\subsection{Spectral Properties of the Sub-Markov Semigroup}
\label{sec:Compact}
To prove the quasi-ergodicity of $(Z_t)_{t\ge0}$ with killing outside of $E$, we study the spectrum of $(Q_t)_{t\ge0}$. 
First, we need the following lemma. 

\begin{lemma}
\label{lemma:Compact}
Let Assumptions \ref{assn:Compact}, \ref{assn:Irreducible}, \&~\ref{assn:tau}~hold. 
Then, we have the following. 
\begin{enumerate}[(i)]
\item $(Q_t)_{t\ge0}$ extends to a strongly continuous semigroup of compact operators on $L^2(E,\mu)$. 
\item $(Q_t)_{t\ge0}$ is strongly Feller. 
\end{enumerate}
\begin{proof}
(i) We first prove strong continuity. 
For $f\in BM(E)$, let $\overline{f}$ be the extension of $f$ to $H$ by zero. 
For arbitrary $p\ge1$, $f\in BM(E)$, and all $t>0$, $x\in E$, we have 
\[
\abs*{Q_tf(x)}^p\,\le\,\ExpOp_x\left[\abs*{f(Z_t)}^p1_{t<\tau}\right]\,\le\,\ExpOp_x\left[\abs*{\overline{f}(Z_t)}^p\right]\,\le\, P_t\left(\abs{\overline{f}}^p\right)(x), 
\]
by Jensen's inequality. 
Hence, by the $\mu$-invariance of $(P_t)_{t\ge0}$, 
\[
\begin{aligned}
\int_E\abs*{Q_tf(x)}^p\,\mu(dx)\,&\le\,\int_H P_t\left(\abs*{\overline{f}}^p\right)(x)\,\mu(dx) \\
&=\,\int_H\abs*{\overline{f}(x)}^p\,\mu(dx)\,=\,\int_E \abs*{f(x)}^p\,\mu(dx). 
\end{aligned}
\]
Thus $\norm*{Q_tf}_{L^p_\mu}\le\norm*{f}_{L^p_\mu}$ for $p\ge1$, so $\norm*{Q_tf}_{L^\infty_\mu}\le\norm*{f}_{L^\infty_\mu}$. 
It follows that $Q_tf$ is bounded $\mu$-almost surely by $f\in BM(E)$ uniformly in $t\ge0$, and we may apply the dominated convergence theorem to find that 
\[
\frac{1}{2}\lim_{t\rightarrow0}\norm*{Q_tf-f}_{L^2_\mu}\,\le\,\lim_{t\rightarrow0}\left(\int_E f(x)^2\,\mu(dx) -\int_E f(x)Q_tf(x)\,\mu(dx)\right)\,=\,0. 
\]
Compactness of $Q_t$ follows from from \cite[Theorem 2.3]{AB80}, observing that $Q_t$ is dominated by $P_t$, in the sense that for all positive $f\in L^p(E,\mu)$, 
\[
Q_tf(x)\,=\,\ExpOp_x\left[f(Z_t)1_{t<\tau}\right]\,\le\,\ExpOp_x\left[f(Z_t)\right]\,=\,P_tf(x),\qquad x\in E,\,\,t>0. 
\]

(ii) The following is similar to an argument found in \cite{C85}. 
Fix arbitrary $t>s>0$ and a bounded measurable function $f:E\rightarrow\R$, and define 
\[
\psi_{t,s}(x)\,\coloneqq\,\ExpOp_x\left[f(X_{t-s}),\,t-s<\tau\right]. 
\]
Then $Q_tf(x)=Q_sQ_{t-s}f(x)=\ExpOp_x\left[\psi_{t,s}(X_s)1_{s<\tau}\right]$, and for all $x\in E$ it holds that 
\begin{equation}
\label{eq:Qeq}
\abs*{Q_tf(x)-P_s\psi_{t,s}(x)}\,=\,\ExpOp_x\left[\abs*{\psi_{t,s}(X_s)}1_{s>\tau}\right]\,\le\,\norm*{f}_{L_\mu^\infty}\pr_x\left[s>\tau\right]. 
\end{equation}
Note that the right hand side of \eqref{eq:Qeq}~tends to zero as $s\rightarrow0$ by continuity of $(Z_t)_{t\ge0}$. 
Since $\psi_{t,s}$ is a bounded measurable function and $(P_t)_{t\ge0}$ is strong Feller, we see that $P_s\psi_{t,s}\in C_b(E)$. 
Hence for any sequence $(x_n)_{n\in\N}\subset E$ such that $x_n\rightarrow x\in E$, for any $\epsilon>0$ we can find small $s>0$ and $N>0$ such that if $n>N$, then 
\[
\begin{aligned}
\abs*{Q_tf(x_n)-Q_tf(x)}\,&\le\,\abs*{Q_tf(x_n)-P_s\psi_{t,s}(x_n)} + \abs*{P_s\psi_{t,s}(x_n)-P_s\psi_{t,s}(x)} \\
&\qquad\qquad + \abs*{P_s\psi_{t,s}(x) - Q_tf(x)} \\
&\le\,\epsilon. 
\end{aligned}
\]
Therefore for arbitrary bounded measurable $f:E\rightarrow\R$ and any $t>0$, we have that $Q_tf:E\rightarrow\R$ is continuous. 
\end{proof}
\end{lemma}

Before proceeding, we must borrow a few ideas from the theory of Banach lattices. 
Since these concepts are not the main focus of this paper, we refer to \cite[Section 10.3]{BFR17}~for the relevant definitions. 
The proof of the following result is similar to an argument in \cite[Theorem 4.5]{CLMR21}, with a subtle difference due to the $L^2$, continuous time setting in which we work. 

\begin{lemma}
\label{lemma:sIrreducible}
Under Assumptions \ref{assn:Compact}~\&~\ref{assn:Irreducible}, $(Q_t)_{t\ge0}$ is an ideal irreducible semigroup on the Banach lattice $L^2(E,\mu)$, in the sense that for each $t>0$, the only closed ideals in $L^2(E,\mu)$ that are $Q_t$-invariant are $L^2(E,\mu)$ and $\{0\}$. 
\begin{proof}
Let $I$ be a closed ideal in $L^2(E,\mu)$. 
By \cite[Proposition 10.15]{BFR17}, there must exist a measurable set $A\subset E$ such that 
\[
I\,=\,\left\{f\in L^2(E,\mu)\,:\,f\rvert_A\,=\,0\,\,\text{$\mu$-almost surely}\right\}. 
\]
If $\mu(A)=\mu(E)$, then $I$ consists solely of the zero function in $L^2(E,\mu)$, while $\mu(A)=0$ implies that $I=L^2(E,\mu)$. 
Suppose that $0<\mu(A)<\mu(E)$, and take $f\in I/\{0\}$ such that $f\ge0$. 
Hence there exists a real number $\epsilon>0$ and a set $B\subset E/A$ of positive $\mu$-measure such that $f(x)>\epsilon$ for $\mu$-almost all $x\in B$. 
If we may take $B$ to be open, then by Assumption \ref{assn:tau}~we observe that for $\mu$-almost all $x\in E$, 
\begin{equation}
\label{eq:OpenIrreducibility}
Q_t1_B(x)\,=\,\ExpOp_x\left[1_B(Z_t)1_{t<\tau}\right]\,=\,\pr_x\left[Z_t\in B,\,t<\tau\right]\,\ge\,\pr_x\left[t<\tau_B\right]\,>0\,. 
\end{equation}
It follows that $Q_tf>\epsilon Q_t1_B$ holds $\mu$-almost surely in $E$, and in particular that $Q_tf$ cannot be contained in $I$. 

If $B$ cannot be taken to be open, let $t>s>0$, and note that $Q_s1_B$ is a continuous function, by Lemma \ref{lemma:Compact}. 
Moreover, by the strong continuity of $(Q_t)_{t\ge0}$, we may take $s>0$ such that $Q_s1_B$ is nonnegative and nonzero. 
Therefore, $Q_s1_B$ has nontrivial open support. 
Hence we may find an open nonempty set $B'\subset E$ and $\epsilon'>0$ such that $Q_s1_B(x)\ge\epsilon'1_{B'}(x)$ for $\mu$-almost all $x\in E$. 
Then, using the positivity of $(Q_t)_{t\ge0}$ we find that for all $t>0$ and $s>0$ sufficiently small, we have 
\[
Q_tf\,\ge\,\epsilon Q_t1_B\,=\,\epsilon Q_{t-s}Q_s1_B\,\ge\,\epsilon\epsilon'Q_{t-s}1_{B'}\,>\,0, 
\]
holds $\mu$-almost surely in $E$, where the last inequality holds as in \eqref{eq:OpenIrreducibility}~with $B$ replaced by $B'$. 
In particular, we again see that $Q_tf(x)>0$ for $\mu$-almost all $x\in A$, and hence $Q_tf$ is not in $I$. 

Therefore if $0<\mu(A)<\mu(E)$, it is impossible for $I$ to be $Q_t$-invariant. 
Hence for arbitrary $t>0$, the only $Q_t$-invariant ideals in $L^2(E,\mu)$ are $\{0\}$ and $L^2(E,\mu)$. 
\end{proof}
\end{lemma}

The generator of $(P_t)_{t\ge0}$ in $L^2(H,\mu)$ is denoted $\mathcal{L}$, while the generator of $(Q_t)_{t\ge0}$ in $L^2(E,\mu)$ is denoted $\mathcal{L}_E$. 
For each $t>0$, $Q_t$ has an adjoint in $L^2(E,\mu)$, denoted $Q_t^*$. 
Define 
\[
s(Q_t)\,\coloneqq\,\sup\spec Q_t,\qquad s(\mathcal{L}_E)\,\coloneqq\,\sup\spec\mathcal{L}_E. 
\]
Given Lemma \ref{lemma:Compact}~\&~\ref{lemma:sIrreducible}, we have the following classical result. 

\begin{proposition}
\label{thm:Simple}
$s(\mathcal{L}_E)$ is a simple eigenvalue of $\mathcal{L}_E$ with positive eigenvector. 
\begin{proof}
Since for each $t>0$ the operator $Q_t$ is compact, $\mathcal{L}_E$ must have compact resolvent \cite[Theorem II.4.29]{ENB00}. 
Since each $Q_t$ is positive, $\mathcal{L}_E$ is resolvent positive \cite[Corollary 11.4]{BFR17}. 
Moreover, since $Q_t$ is a positive compact irreducible operator, \cite{DP86}~implies that $s(Q_t)>0$, and therefore by the spectral mapping theorem \cite[Theorem 2.4]{P12} 
\[
s(\mathcal{L}_E)\,=\,\frac{1}{t}\ln s(Q_t)\,>\,-\infty. 
\] 
By \cite[Proposition 12.15]{BFR17}, there exists $\varphi\in D(A)$ such that $\varphi>0$ and $\mathcal{L}_E\varphi=s(\mathcal{L}_E)\varphi$. 
As $s(\mathcal{L}_E)$ is an eigenvalue of $\mathcal{L}_E$, it must be a pole of the resolvent operator of $\mathcal{L}_E$. 
By \cite[Theorem 14.12(d)]{BFR17}, the order of this pole has multiplicity equal to one. 
Consequently, the eigenspace of $s(\mathcal{L}_E)$ is equal to $\text{span}\{\varphi\}$. 
\end{proof}
\end{proposition}

\subsection{Quasi-Ergodicity of Markov Processes on Separable Hilbert Space}
\label{sec:QuasiQ}
Proposition \ref{thm:Simple}~in hand, we are able to prove results on the quasi-ergodicity of $(Z_t)_{t\ge0}$. 
Let $\varphi,\,\varphi^*$, be the eigenfunctions from Proposition \ref{thm:Simple}, normalized such that 
\[
\int_E\varphi(x)\,\mu(dx)\,=\,\int_E\varphi^*(x)\,\mu(dx)\,=\,1, 
\]
and define 
\[
M\,\coloneqq\,\inner*{\varphi}{\varphi^*}_{L_\mu^2}\,\in\,(0,1]. 
\]
We have the following lemma on the asymptotic behaviour of $(Q_t)_{t\ge0}$. 
The proof idea is similar to \cite[Lemma 3]{B57}. 
In the case where $H$ is finite dimensional and $Q_t$ possesses a uniformly bounded integral kernel with respect to $\mu$, similiar results are proven in \cite[Theorem 3]{P85}, and also in \cite{KS08,ZLS14}. 
Our weaker assumptions lead us to convergence results in $L^2$ and almost surely, rather than uniformly. 

\begin{lemma}
\label{lemma:}
Let Assumptions \ref{assn:Compact}, \ref{assn:Irreducible}, \&~\ref{assn:tau}~hold. 
Then, there exists $\gamma>0$ such that for all $f\in L^2(E,\mu)$ we have 
\begin{equation}
\label{eq:AsymptoticQt}
\norm*{e^{\lambda_1t}Q_tf(\cdot) - M^{-1}\int\varphi(\cdot)\varphi^*(y)f(y)\,\mu(dy)}_{L^2_\mu}\,\le\,e^{-\gamma t}\left(1+M^{-1}\right)\norm*{f}_{L^2_\mu}. 
\end{equation}
Moreover, $e^{\lambda_1t}Q_tf(\cdot)$ converges $\mu$-almost surely to $M^{-1}\int\varphi(\cdot)\varphi^*(y)f(y)\,\mu(dy)$. 
\begin{proof}
By Proposition \ref{thm:Simple}~we may split $L^2(E,\mu)=\mathcal{H}_1\otimes\mathcal{H}_{rem}$, where $\mathcal{H}_1=\text{span}\{\varphi\}$ and both of $\mathcal{H}_1,\,\mathcal{H}_{rem}$ are $Q_t$-invariant, see Deimling \cite[Theorem 8.9]{D10}. 
Since $e^{-\lambda_1t}\,=\,\sup\spec(Q_t)$ and there is a gap between $e^{-\lambda_1t}$ and the rest of the spectrum of $Q_t$, there exists $\gamma>0$ such that 
\begin{equation}
\label{eq:Gap}
e^{\lambda_1t}\norm*{Q_t\rvert_{\mathcal{H}_{rem}}}_{L^2_\mu}\,\le\,e^{-\gamma t}\quad\text{ for large }\,\,t>0. 
\end{equation}
Moreover, for $f\in L^2(E,\mu)$ there exists $\psi_f\in\mathcal{H}_{rem}$ and $c_f\in\R$ such that 
\[
f\,=\,c_f\varphi + \psi_f. 
\]
Observing that 
\[
\begin{aligned}
0\,&=\,\lim_{t\rightarrow\infty}\int_E \left(e^{\lambda_1t}Q_tf(x) - c_f\varphi(x)\right)\varphi^*(x) \,\mu(dx) \\
&=\,\lim_{t\rightarrow\infty}\int_E e^{\lambda_1t}f(x)Q_t^*\varphi^*(x)\,\mu(dx) - c_fM \\
&=\,\lim_{t\rightarrow\infty}\int_E f(x)\varphi^*(x)\,\mu(dx) -c_fM, 
\end{aligned}
\]
we obtain $c_f=M^{-1}\inner*{f}{\varphi^*}_{L_\mu^2}$. 
By H{\"o}lder's inequality, 
\[
\abs*{c_f}\,\le\,M^{-1}\norm*{f}_{L^2_\mu},\qquad\norm*{\psi_f}_{L^2_\mu}\,\le\,\left(1+M^{-1}\right)\norm*{f}_{L^2_\mu}, 
\]
and hence 
\[
\begin{aligned}
&\norm*{e^{\lambda_1t}Q_tf(x) - M^{-1}\int_E\varphi(x)\varphi^*(y)f(y)\,\mu(dy)}_{L^2_\mu}\\
&\qquad=\,\norm*{e^{\lambda_1t}c_fQ_t\varphi(x) + e^{\lambda_1t}Q_t\psi_f - M^{-1}\varphi(x)\inner*{f}{\varphi^*}_{L^2_\mu}}_{L^2_\mu}\\
&\qquad=\,\norm*{\left(c_f- M^{-1}\inner*{f}{\varphi^*}_{L^2_\mu}\right)\varphi(x) + e^{\lambda_1t}Q_t\psi_f}_{L^2_\mu} \\
&\qquad=\,\norm*{e^{\lambda_1t}Q_t\psi_f}_{L^2_\mu}\,\le\,e^{-\gamma t}\left(1+M^{-1}\right)\norm*{f}_{L^2_\mu}, 
\end{aligned}
\]
completing the proof of \eqref{eq:AsymptoticQt}. 

To prove that $e^{\lambda_1t}Q_tf(\cdot)$ converges $\mu$-almost surely to $M^{-1}\int\varphi(\cdot)\varphi^*(y)f(y)\,\mu(dy)$, fix $f\in L^2(E,\mu)$ and define 
\[
A_t(\cdot)\,\coloneqq\,\abs*{e^{\lambda_1t}Q_tf(\cdot)-M^{-1}\int\varphi(\cdot)\varphi^*(y)f(y)\,\mu(dy)}, 
\]
so that $\norm*{A_t}_{L^2_\mu}\,\le\, C_fe^{-\gamma t}$ for some $C_f>0$. 
By Chebyshev's inequality, for any $n\in\N$ and $t,s>0$ we have 
\[
\begin{aligned}
\mu\left\{x\in E\,:\,\abs*{A_t(x)-A_s(x)}>n^{-1}\right\}\,&\le\,n^2\norm*{A_t-A_s}^2_{L^2_\mu}\\
&\le\, 2n^2\min\left\{\norm*{A_s}_{L_\mu^2},\norm*{A_t}_{L^2_\mu}\right\} \\
&\le\, 2n^2C_f\exp\left(-\gamma\min\left\{s,t\right\}\right). 
\end{aligned}
\]
Therefore $(A_t)_{t\ge0}$ is Cauchy $\mu$-almost surely, and so converges $\mu$-almost surely. 
By \eqref{eq:AsymptoticQt}, the $\mu$-almost sure limit of $(A_t)_{t\ge0}$ must be zero, completing the proof. 
\end{proof}
\end{lemma}

We are now ready to prove the main result of this section. 
The following results consist of ``conditonal'' limiting theorems for 1.~the distribution of $(Z_t)_{t\ge0}$, 2.~the time-averaged dynamics of $(Z_t)_{t\ge0}$, and 3.~the distribution of a sequence of measurements of $(Z_t)_{t\ge0}$ taken at different times prior to killing. 
Compare the following results with those of \cite[Section 3]{ZLS14}, noting in particular that we work in an $L^p$ framework due to the fact that we make no assumptions on the existence or boundedness of a density for the transition kernel of $(Q_t)_{t\ge0}$. 

\begin{theorem}
\label{thm:quasiErgodic} 
Let Assumptions \ref{assn:Compact}, \ref{assn:Irreducible}, \&~\ref{assn:tau}~hold. 
Defining $\alpha(dx)\coloneqq\varphi^*(x)\mu(dx)$ and $\beta(dx)\coloneqq\varphi(x)\varphi^*(x)\mu(dx)$, the following results hold. 
\begin{enumerate}
\item For $f\in L^2(E,\mu)$ we have 
\begin{equation}
\label{eq:QSD}
\lim_{t\rightarrow\infty}\ExpOp_x\left[f(Z_t)\,|\,t<\tau\right]\,=\,\alpha(f) 
\end{equation}
in $L^2_\mu$ and $\mu$-almost surely. 
Moreover, $\alpha$ is the unique quasi-stationary measure of $(Z_t)_{t\ge0}$. 
That is, $\alpha$ is the only measure on $E$ such that for any measurable $A\subset E$, 
\[
\pr_\alpha\left[Z_t\in A\,|\,t<\tau\right]\,=\,\alpha(A),\qquad t\,>\,0. 
\]
\item For arbitrary $f\in L^2(E,\mu)$, $\epsilon>0$, and $\mu$-almost all $x\in E$, we have 
\begin{equation}
\label{eq:QED}
\lim_{t\rightarrow\infty}\pr_x\left[\abs*{\frac{1}{t}\int_0^tf(Z_s)\,ds - \beta(f)}>\epsilon\,\big|\,t<\tau\right]\,=\,0. 
\end{equation}
Moreover, $\beta$ is the unique quasi-ergodic measure of $(Z_t)_{t\ge0}$ on $E$.
\item For $f,g\in L^2(E,\mu)$ and $0<a<b<1$ we have  in $L^2_\mu$ and $\mu$-almost surely that 
\begin{align}
\lim_{t\rightarrow\infty}\ExpOp_x\left[f(Z_{at})g(Z_t)\,|\,t<\tau\right]\,&=\,\beta(f)\alpha(g), 
\label{eq:DoubleQSD}\\
\lim_{t\rightarrow\infty}\ExpOp_x\left[f(Z_{at})g(Z_{bt})\,|\,t<\tau\right]\,&=\,\beta(f)\beta(g)
\label{eq:DoubleQSD2}. 
\end{align}
\item For $f\in L^2(E,\mu)$, it holds in $L^2_\mu$ and $\mu$-almost surely that 
\begin{equation}
\label{eq:DoubleQuasiErgodic}
\lim_{t\rightarrow\infty}\lim_{T\rightarrow\infty}\ExpOp_x\left[f(Z_t)\,|\,T<\tau\right]\,=\,\beta(f). 
\end{equation}
\end{enumerate}
\begin{proof}
To see \eqref{eq:QSD}, note that Lemma \ref{lemma:}~implies 
\[
\begin{aligned}
\lim_{t\rightarrow\infty}\ExpOp_x\left[f(Z_t)\,|\,t<\tau\right]\,&=\,\lim_{t\rightarrow\infty}\frac{Q_tf(x)}{Q_t1(x)}\\
&=\,\frac{M^{-1}\int_E\varphi(x)\varphi^*(y)f(y)\,\mu(dy)}{M^{-1}\int_E\varphi(x)\varphi^*(y)\,\mu(dy)}\,=\,\alpha(f), 
\end{aligned}
\]
the above limits holding $\mu$-almost surely and in $L^2_\mu$. 
To see that $\alpha$ is in fact a quasi-stationary measure, we compute for $\mu$-almost all $x\in E$ 
\[
\begin{aligned}
\ExpOp_\alpha\left[f(Z_t)\,|\,t<\tau\right]\,&=\,\frac{\int_E Q_tf(y)\,\alpha(dy)}{\pr_{\alpha}[t<\tau]}\\
&=\, \frac{1}{\pr_{\alpha}[t<\tau]}\lim_{s\rightarrow\infty}\frac{Q_s(Q_tf)(x)}{\pr_x[s<\tau]} \\ 
&=\, \frac{1}{\pr_{\alpha}[t<\tau]}\lim_{s\rightarrow\infty}\frac{Q_s(Q_tf)(x)}{\pr_x[s+t<\tau]}\frac{\pr_x[s+t<\tau]}{\pr_x[s<\tau]}\\
&=\, \frac{1}{\pr_{\alpha}[t<\tau]}\lim_{s\rightarrow\infty}\ExpOp_x\left[f(X_{s+t})\,|s+t<\tau\right]\ExpOp_x[Q_t1(X_s)\,|\,s<\tau] \\ 
&=\,\frac{1}{\pr_{\alpha}[t<\tau]}\int_E f(y)\,\alpha(dy)\int_E Q_t1(y)\,\alpha(dy) \,=\,\int_Ef(y)\,\alpha(dy). 
\end{aligned}
\]
Supposing $\alpha_0$ were a second quasi-stationary measure, then 
\[
\int_Ef(y)\,\alpha_0(dy)\,=\,\lim_{t\rightarrow\infty}\ExpOp_{\alpha_0}\left[f(Z_t)\,|\,t<\tau\right]\,=\,\int_Ef(y)\,\alpha(dy), 
\]
so that $\alpha_0=\alpha$ by duality. 

We now prove \eqref{eq:DoubleQSD}.  
First note that, by the Markov property, 
\begin{equation}
\label{eq:QQMarkov}
\ExpOp_x\left[f(Z_{at})g(Z_t)\,|\,t<\tau\right]\,=\,\frac{e^{\lambda_1t}Q_{at}\left(f(\cdot)Q_{t-at}g(\cdot)\right)(x)}{e^{\lambda_1t}Q_t1(x)}. 
\end{equation}
From Lemma \ref{lemma:}, we have that $\left(e^{\lambda_1t}Q_t1(x)\right)$ converges to $M^{-1}\varphi(x)$ as $t\rightarrow\infty$  $\mu$-almost surely and in $L^2_\mu$. 
We now show that the limit of the numerator in \eqref{eq:QQMarkov}~is $M^{-1}\varphi(x)\beta(f)\alpha(g)$ in $L^2_\mu$ and $\mu$-almost surely. 
To see this, we define 
\[
h_t(x)\,\coloneqq f(x) e^{\lambda_1(t-at)}Q_{t-at}g(x), 
\]
and compute 
\begin{equation}
\begin{aligned}
&e^{\lambda_1t}Q_{at}(f(\cdot)Q_{t-at}g(\cdot))(x) - M^{-1}\varphi(x)\beta(f)\alpha(g)\\
&\qquad=\,
e^{\lambda_1t}Q_{at}(f(\cdot)Q_{t-at}g(\cdot))(x) - \varphi(x)M^{-2}\int\varphi(y)\varphi^*(y)f(y)\,\mu(dy)\int\varphi^*(z)g(z)\,\mu(dz) \\
&\qquad=\, e^{\lambda_1at}Q_{at}h_t(x) - M^{-1}\int\varphi(x)\varphi^*(y)f(y) M^{-1}\int\varphi(y)\varphi^*(z)g(z)\,\mu(dz)\,\mu(dy) \\
&\qquad=\, e^{\lambda_1at}Q_{at}h_t(x) \\
&\qquad\qquad- M^{-1}\int\varphi(x)\varphi^*(y)f(y)\bigg(M^{-1}\int\varphi(y)\varphi^*(z)g(z) \\
&\qquad\qquad\qquad\qquad\qquad\qquad\qquad\qquad\qquad - e^{\lambda_1(t-at)}Q_{t-at}g(y) + e^{\lambda_1(t-at)}Q_{t-at}g(y)\bigg)\,\mu(dy) \\
&\qquad=\, e^{\lambda_1at}Q_{at}h_t(x) - M^{-1}\int\varphi(x)\varphi^*(y)h_t(y)\,\mu(dy) \\
&\qquad\qquad+ M^{-1}\int\varphi(x)\varphi^*(y)f(y)\left( e^{\lambda_1(t-at)}Q_{t-at}g(y) - M^{-1}\int\varphi(y)\varphi^*(z)g(z)\,\mu(dz) \right)\,\mu(dy). 
\end{aligned}
\end{equation}
Taking the $L^2_\mu$ norm and applying H{\"o}lder's inequality, we have 
\begin{equation}
\label{eq:Lim}
\begin{aligned}
&\norm*{e^{\lambda_1at}Q_{at}h_t(x) - M^{-1}\int\varphi(x)\varphi^*(y)h_t(y)\,\mu(dy)}_{L^2_\mu} \\
&\quad+ M^{-1}\norm*{\varphi}_{L^2_\mu}\abs*{\int\varphi^*(y)f(y)\left( e^{\lambda_1(t-at)}Q_{t-at}g(y) - M^{-1}\int\varphi(y)\varphi^*(z)g(z)\,\mu(dz) \right)\,\mu(dy)} \\
&\quad\le\, 
\norm*{e^{\lambda_1at}Q_{at}h_t(x) - M^{-1}\int\varphi(x)\varphi^*(y)h_t(y)\,\mu(dy)}_{L^2_\mu} \\
&\quad\qquad+ M^{-1}\norm*{\varphi}_{L^2_\mu}\norm*{\varphi^*f}_{L^2_\mu}\norm*{e^{\lambda_1(t-at)}Q_{t-at}g(\cdot) - M^{-1}\int\varphi(\cdot)\varphi^*(z)g(z)\,\mu(dz)}_{L^2_\mu} \\
&\qquad\le\, e^{-\gamma t}(1+M^{-1})\norm*{h_t}_{L^2_\mu} + M^{-1}\norm*{\varphi}_{L^2_\mu}\norm*{\varphi^*f}_{L^2_\mu}e^{-\gamma(t-at)}(1+M^{-1})\norm*{g}_{L^2_\mu}. 
\end{aligned}
\end{equation}
Note that 
\[
\begin{aligned}
\norm*{h_t}_{L^2_\mu}\,&=\,\norm*{f(\cdot)e^{\lambda_1(t-at)}Q_{t-at}g(\cdot)}_{L^2_\mu} \\
&\le\, \norm*{f}_{L^2_\mu}e^{\lambda_1(t-at)}\norm*{Q_{t-at}g}_{L^2_\mu}\,\le\,\norm*{f}_{L^2_\mu}\norm*{g}_{L^2_\mu}, 
\end{aligned}
\]
so that \eqref{eq:Lim}~tends to zero as $t\rightarrow\infty$. 
Hence  
\begin{equation}
\label{eq:L2Convergence}
\lim_{t\rightarrow\infty} e^{\lambda_1t}Q_{at}(f(\cdot)Q_{t-at}g(\cdot))(x)\,=\,M^{-1}\varphi(x)\beta(f)\alpha(g).  
\end{equation}
in $L^2_\mu$. 
Since the rate of convergence in \eqref{eq:L2Convergence}~is exponential, the same argument as in Lemma \ref{lemma:}~can be used to show that the convergence in \eqref{eq:L2Convergence}~holds $\mu$-almost surely. 
We have therefore proven that $\lim_{t\rightarrow\infty}\ExpOp_x\left[f(Z_{at})g(Z_t)\,|\,t<\tau\right]\,=\,\beta(f)\alpha(g)$ in $L^2_\mu$ and $\mu$-almost surely. 
The proofs of \eqref{eq:DoubleQSD2}~\&~\eqref{eq:DoubleQuasiErgodic}~are similar. 

To prove \eqref{eq:QED}, note that from \eqref{eq:QQMarkov}~we have 
\[
\ExpOp_x\left[\frac{1}{t}\int_0^tf(Z_s)\,ds\,|\,t<\tau\right]\,=\,\frac{\frac{1}{t}\int_0^t e^{\lambda_1t}Q_s\left(f(\cdot)Q_{t-s}1(\cdot)\right)(x)\,ds}{e^{\lambda_1t}Q_t1(x)}. 
\]
Then, defining $h_t(x)\coloneqq f(x)Q_{t-s}1(x)$, observe that 
\[
\begin{aligned}
&\frac{1}{t}\int_0^t e^{\lambda_1t}Q_s\left(f(\cdot)Q_{t-s}1(\cdot)\right)(x)\,ds - M^{-1}\varphi(x)\beta(f)\\
&\qquad=\,\frac{1}{t}\int_0^t e^{\lambda_1s}Q_s\left(f(\cdot)e^{\lambda_1(t-s)}Q_{t-s}1(\cdot)\right)(x) - M^{-1}\varphi(x)\beta(f)\,ds\\
&\qquad=\,\frac{1}{t}\int_0^t e^{\lambda_1s}Q_s\left(f(\cdot)e^{\lambda_1(t-s)}Q_{t-s}1(\cdot)\right)(x) \\
&\qquad\qquad
- M^{-1}\int\varphi(x)\varphi^*(y)f(y)\, M^{-1}\int\varphi(y)\varphi^*(z)1(z)\,\mu(z)\mu(dy) 
\,ds\\
&\qquad=\, \frac{1}{t}\int_0^te^{\lambda_1s}Q_sh_t(x) - M^{-1}\int\varphi(x)\varphi^*(y)h_t(y)\,\mu(dy) \\ 
&\qquad\qquad - M^{-1}\int\varphi(x)\varphi^*(y)\left(M^{-1}\int\varphi(y)\varphi^*(z)1(z)\mu(dz) - e^{\lambda_1(t-s)}Q_{t-s}1(y)\right)\,ds. 
\end{aligned}
\]
As the $L^2_\mu$ norm of $h_t$ is bounded uniformly in $t>0$, the same argument as used to prove \eqref{eq:DoubleQSD}~implies that 
\[
\lim_{t\rightarrow\infty}\norm*{\ExpOp_x\left[\frac{1}{t}\int_0^tf(Z_s)\,ds\,|\,t<\tau\right]-\beta(f)}_{L^2_\mu}\,=\,0. 
\]
Since the above convergence can be seen to occur at an exponential rate, 
\[
\lim_{t\rightarrow\infty}\ExpOp_x\left[\frac{1}{t}\int_0^tf(Z_s)\,ds\,|\,t<\tau\right]\,=\,\beta(f) 
\]
$\mu$-almost surely, again using the argument as in Lemma \ref{lemma:}. 
Applying Markov's inequality completes the proof. 
\end{proof}
\end{theorem}

Now, we prove that the rate of convergence to the quasi-stationary measure $\alpha$ in Theorem \ref{thm:quasiErgodic}.1~is exponential in $L^2(E,\alpha^*)$, where $\alpha^*(dx)\coloneqq\varphi(x)\mu(dx)$. 
If $(Z_t)_{t\ge0}$ is reversible with respect to $\mu$, we of course have $\alpha=\alpha^*$. 

\begin{theorem}
\label{thm:expQSD}
Let Assumptions \ref{assn:Compact}, \ref{assn:Irreducible}, \&~\ref{assn:tau}~hold, and define the probability measure $\alpha^*(dx)\coloneqq\varphi(x)\mu(dx)$. 
Then, for any $f\in L^2(E,\mu)$ there exist $T>0$ and $C_f>0$ such that 
\[
\norm*{\ExpOp_x\left[f(Z_t)\,|\,t<\tau\right]-\alpha(f)}_{L^1_{\alpha^*}}\,\le\,C_fe^{-\gamma t}\qquad\forall t\ge T. 
\]
\begin{proof}
First, note that since $\frac{Q_tf}{Q_t1}$ converges to $\alpha(f)$ in $L^2_\mu$, for any $K>0$ there exists $T>0$ such that for all $t\ge T$, we have 
\[
\begin{aligned}
\norm*{\frac{M^{-1}}{Q_t1}Q_tf - e^{\lambda t}Q_tf}_{L_\mu^1}\,&=\,\norm*{\left(M^{-1}\varphi - e^{\lambda t}Q_t1\right)\frac{Q_tf}{Q_t1}}_{L^1_\mu} \\
&\le\, \norm*{e^{\lambda t}Q_t1-M^{-1}\varphi}_{L^2_\mu}\norm*{\frac{Q_tf}{Q_t1}}_{L^2_\mu}\\
&\le\,C_f'e^{-(\gamma-\lambda)t} 
\end{aligned}
\]
for some $C'_f>0$. 
From this we compute 
\[
\begin{aligned}
\norm*{\ExpOp_x\left[f(Z_t)\,|\,t<\tau\right]-\alpha(f)}_{L^1_{\alpha^*}}\,&=\,\norm*{\frac{\varphi}{Q_t1}Q_tf - \alpha(f)\varphi}_{L^1_\mu}\\
&\le\,M\norm*{\frac{M^{-1}\varphi}{Q_t1}Q_tf - e^{\lambda t}Q_tf + e^{\lambda t}Q_tf - M^{-1}\alpha(f)\varphi}_{L^1_\mu}\\
&\le\,M\norm*{\frac{M^{-1}\varphi}{Q_t1}Q_tf - e^{\lambda t}Q_tf}_{L^1_\mu} + M\norm*{e^{\lambda t}Q_tf - M^{-1}\alpha(f)\varphi}_{L^2_\mu} \\
&\le\, MC_f'e^{-(\gamma-\lambda)t} + MC''_f e^{-\gamma t}, 
\end{aligned}
\]
where $C_f''$ is as in Lemma \ref{lemma:}, completing the proof. 
\end{proof}
\end{theorem}

Since we have not proven that the rate of convergence to $\alpha$ is uniform, the proof of existence of a $Q$-process found in \cite{CV16}~does not immediately translate to our setting. 
Nevertheless, we prove the existence and uniqueness of the $Q$-process in Theorem \ref{thm:Qproc}~below using arguments inspired by their work. 

\begin{theorem}
\label{thm:Qproc}
Under Assumptions \ref{assn:Compact}, \ref{assn:Irreducible}, \&~\ref{assn:tau}, for $\mu$-almost all $x\in E$ there exists a probability measure $\mathbb{Q}_x$ defined set-wise as the limit 
\begin{equation}
\label{eq:Qx}
\mathbb{Q}_x(A)\,\coloneqq\,\lim_{t\rightarrow\infty}\pr_x(A\,|\,t<\tau),\qquad A\in\mathcal{F}_s,\,\,s\ge0. 
\end{equation}
Moreover, the limit in \eqref{eq:Qx}~holds in $L^2_\mu$. 
With respect to $(\mathbb{Q}_x)_{x\in E}$, the solution process $(Z_t)_{t\ge0}$ is a homogeneous Markov process in $E$ with unique ergodic measure $\beta$ 
\begin{proof}
For $x\in E$ and $t>0$, we define an auxiliary probability measure $\mathbb{Q}_x^t$ on $\Omega$ as 
\[
\mathbb{Q}_x^t(d\omega)\,\coloneqq\,\frac{1_{t<\tau}(\omega)}{\ExpOp_x[1_{t<\tau}]}\pr_x(d\omega). 
\]
By the Markov property of $(Z_t)_{t\ge0}$ with respect to $\pr_x$, for fixed $s>0$ and any $t\ge s$ 
\begin{equation}
\label{eq:fracExpOp}
\frac{\ExpOp_x\left[1_{t<\tau}\,|\,\mathcal{F}_s\right]}{\ExpOp_x\left[1_{t<\tau}\right]}\,=\,\frac{1_{s<\tau}\pr_{Z_s}\left[t-s<\tau\right]}{\pr_x[t<\tau]}.
\end{equation}
By Lemma \ref{lemma:}~with $\pr_x[t<\tau]=Q_t1(x)$, we have in $L^2_\mu$ and $\mu$-almost surely that 
\[
M^{-1}\varphi(x)\,=\,\lim_{t\rightarrow\infty}\frac{\pr_x[t<\tau]}{\pr_\alpha[t<\tau]}. 
\]
Hence for fixed $s>0$, \eqref{eq:fracExpOp}~converges $\mu$-almost surely as $t\rightarrow\infty$ to 
\[
M_s(x)\,\coloneqq\,1_{s<\tau}e^{\lambda_1s}\frac{\varphi(Z_s)}{\varphi(x)}. 
\]
Also, we may compute 
\begin{equation}
\label{eq:M1} 
\ExpOp_x[M_s(x)]\,=\,e^{\lambda_1s}\varphi(x)^{-1}\ExpOp_x\left[1_{s<\tau}\varphi(Z_s)\right]\,=\,e^{\lambda_1s}\varphi(x)^{-1}Q_s\varphi(x)\,=\,1. 
\end{equation}

Now, we claim that the definition of $(M_t)_{t\ge0}$ and \eqref{eq:M1}~imply that for each $A_s\in\mathcal{F}_s$, $s\ge0$, it holds that 
\begin{equation}
\label{eq:claimM}
\lim_{t\rightarrow\infty}\mathbb{Q}_x^t(A_s)\,=\,\ExpOp_x\left[1_{A_s}M_s\right]\qquad\text{ $\mu$-almost surely}. 
\end{equation}
This claim is similar to \cite[Theorem 2.1]{RVY06}. 
For $s\ge0$, $A_s\in\mathcal{F}_s$, and $\omega\in A_s$ note that $1_{t<\tau}(\omega)=\ExpOp_x[1_{t<\tau}|\mathcal{F}_s](\omega)$, and hence 
\[
\begin{aligned}
\lim_{t\rightarrow\infty}\int_{A_s}\norm*{\frac{1_{t<\tau}(\omega)}{\ExpOp_x[1_{t<\tau}]}-M_s(\omega)}\,\pr_x(d\omega)\,&\le\, \lim_{t\rightarrow\infty}\int_{A_s}\norm*{\frac{1_{t<\tau}(\omega)}{\ExpOp_x[1_{t<\tau}]}}\,\pr_x(d\omega) \\
&\qquad\qquad- \lim_{t\rightarrow\infty}\int_{A_s}\norm*{\frac{\ExpOp[1_{t<\tau}|\mathcal{F}_s](\omega)}{\ExpOp[1_{t<\tau}]}}\,\pr_x(d\omega) \\
&=\,0 
\end{aligned}
\]
for $\mu$-almost all $x\in E$. 
By Scheff{\'e}'s Lemma \cite{S47}, it then follows that 
\[
\lim_{t\rightarrow\infty}\mathbb{Q}_x^t(A_s)\,=\,\int_{A_s}M_s(\omega)\,\pr_x(d\omega)\,=\,\ExpOp_x[1_{A_s}M_s] \qquad\text{ $\mu$-almost surely}, 
\]
proving \eqref{eq:claimM}. 
It follows that $\mathbb{Q}_x$ is well-defined for $\mu$-almost all $x\in E$, and  
\[
\frac{d\mathbb{Q}_x}{d\pr_x}\bigg\rvert_{\mathcal{F}_s}\,=\,M_s(x)\qquad\text{ for $s>0$}. 
\]

We now show that $(Z_t)_{t\ge0}$ is a Markov process with respect to $(\mathbb{Q}_x)_{x\in E}$. 
Indeed, using the definition of conditional expectations we have 
\[
\begin{aligned}
M_s\ExpOp_x^{\mathbb{Q}}\left[f(Z_t)\,|\,\mathcal{F}_s\right]\,&=\,
\ExpOp_x^{\mathbb{Q}}\left[M_tf(Z_t)\,|\,\mathcal{F}_s\right]\\
&=\,\ExpOp_x^{\mathbb{Q}}\left[M_tf(Z_t)\,|\,Z_s\right]\,=\,M_s\ExpOp_x^{\mathbb{Q}}\left[f(Z_t)\,|\,Z_s\right], 
\end{aligned}
\]
the second equality following from the definition of $M_t$ and the Markov property of $(\pr_x)_{x\in E}$. 
This proves that $(Z_t)_{t\ge0}$ is a Markov process with respect to $(\mathbb{Q}_x)_{x\in E}$. 
The fact that $(Z_t)_{t\ge0}$ with respect to $(\mathbb{Q}_x)_{x\in E}$ has a unique ergodic measure given by $\beta$ follows from the fourth statement of Theorem \ref{thm:quasiErgodic} . 
\end{proof}
\end{theorem}

\section{Application to a class of SPDEs}
\label{sec:Examples} 
We now demonstrate that Assumption \ref{assn:Supp}~guarantees that a cutoff version of \eqref{eq:SPDE}~satisfies Assumptions \ref{assn:Compact}, \ref{assn:Irreducible}, \&~\ref{assn:tau}, and therefore possesses unique quasi-stationary and quasi-ergodic measures in $\Gamma_\delta$. 
We then show, in Proposition \ref{prop:CutoffEquivalence}, that the quasi-stationary and quasi-ergodic measure of this cutoff version of \eqref{eq:SPDE}~are a quasi-stationary and quasi-ergodic measure of \eqref{eq:SPDE}~itself. 

Before proceeding, we make a few remarks on the setup for our SPDE and the notation used in this section. 
Since $\Gamma_\delta$ is not necessarily open in $H$, we let $\Gamma_\delta^H$ be the smallest open set in $H$ such that $\Gamma_\delta^H\cap D(N)=\Gamma_\delta$. 
Throughout this section, if $E_1,\,E_2,$ are Banach spaces, the Fr{\'e}chet derivative of a function $f:E_1\rightarrow E_2$ at $x\in E_1$ is denoted $Df(x)$. 

Now, letting $N_\delta:D(N)\rightarrow D(N)$ be a twice continuously Fr{\'e}chet differentiable function such that 
\begin{enumerate}[(i)]
\item $N_\delta(x)=N(x)$ for $x\in\Gamma_\delta$, 
\item for some $\delta'>\delta$ we have $N_\delta(x)=0$ for $x\in D(N)/\Gamma_{\delta'}$, and 
\item $\lip N_\delta=\lip N\rvert_{\Gamma_\delta}=\kappa$, 
\end{enumerate}
we consider 
\begin{equation}
\label{eq:cSPDE}
dX'\,=\,\left(LX'+N_\delta(X')\right)\,dt + \sigma B\,dW. 
\end{equation}
We can see that the exit time of $(X_t')_{t\ge0}$ from $\Gamma_\delta^H$ is equal to the exit time of $(X_t')_{t\ge0}$ from $\Gamma_\delta$, since $\Gamma_\delta^H\cap D(N)=\Gamma_\delta$ and -- as will be seen in Proposition \ref{thm:Solution}~below -- $X_t'\in D(N)$ for all $t>0$. 
Moreover, as $N$ is only defined on $D(N)$, we have $N\rvert_{\Gamma_\delta}=N\rvert_{\Gamma_\delta^H}$. 

Before proceeding, we establish a mild solution theory for \eqref{eq:SPDE}~\&~\eqref{eq:cSPDE}. 
A solution theory of the associated Ornstein-Uhlenbeck process, 
\begin{equation}
\label{eq:OU} 
dY_t\,=\,LY_t\,dt + \sigma B\,dW_t, 
\end{equation}
and the first variational equation of \eqref{eq:cSPDE}~with initial condition $y\in H$, 
\begin{equation}
\label{eq:Variational}
\partial_t\eta_t[y]\,=\,L\eta_t[y] + DN_\delta(X'_t)\eta_t[y],\qquad\eta_0[y]\,=\,y.   
\end{equation}
is also needed. 
The solution theory of \eqref{eq:SPDE}~and \eqref{eq:cSPDE}~follows from a relatively standard fixed point argument, similar to that in \cite[Theorem 7.7]{DPZ14}. 

\begin{proposition}
\label{thm:Solution}
Under Assumptions \ref{assn:PDE}~\&~\ref{assn:SPDE}, we have the following results. 
\begin{enumerate}[(i)]
\item For each $x\in H$, there exists a unique mild solution $(X'_t)_{t\ge0}\subset D(N)$ to \eqref{eq:cSPDE}, 
\begin{equation}
\label{eq:cMild}
X'_t\,=\,\Lambda_tx + \int_0^t\Lambda_{t-s}N_\delta(X'_s)\,ds + \sigma\int_0^t\Lambda_{t-s}B\,dW_s,\qquad t\ge0. 
\end{equation}
\item For each $x\in H$ there exists a unique mild solution $(X_t)_{t\in[0,\tau)}\subset D(N)$ to \eqref{eq:SPDE}, 
\begin{equation}
\label{eq:Mild} 
X_t\,=\,\Lambda_tx + \int_0^t\Lambda_{t-s}N(X_s)\,ds + \sigma\int_0^t\Lambda_{t-s}B\,dW_s,\qquad t\in[0,\tau). 
\end{equation}
\item There exists a unique mild solution to \eqref{eq:OU}~which is continuous in $D(N)$, 
\begin{equation}
\label{eq:Yt}
Y_t\,=\,\Lambda_tY_0 + \sigma\int_0^t\Lambda_{t-s}B\,dW_s,\qquad t\ge0. 
\end{equation}
\item There exists a unique solution $(\eta_t)_{t\ge0}\subset D(N)$ to the first variational equation \eqref{eq:Variational}. 
Letting $X_t'(x)$ denote the solution to \eqref{eq:cSPDE}~with initial condition $x\in H$, $x\mapsto X_t'(x)$ is almost surely continuously Fr{\'e}chet differentiable with $DX_t'(x)y=\eta_t[y]$. 
\end{enumerate}
\begin{proof}
That $(Y_t)_{t\ge0}$ in \eqref{eq:Yt}~is a unique solution to \eqref{eq:OU}~and continuous in $D(N)$ is a consequence of Assumption \ref{assn:SPDE}. 

We now prove existence and uniqueness of mild solutions to \eqref{eq:cSPDE}.  
Since these arguments are largely known under Assumption \ref{assn:SPDE}, we sketch the proof. 
In this case, $N_\delta$ is globally Lipschitz on $D(N)$ with Lipschitz constant $\kappa>0$. 
Fix an arbitrary $T>0$ and $X_0\in D(N)$. 
Let $(Y_t)_{t\ge0}$ be defined as in \eqref{eq:Yt}~with $Y_0=0$, and for some $D(N)$-valued Markov process $(X'_t)_{t\ge0}$ define $V_t\coloneqq X'_t-Y_t$. 
Note that $(X'_t)_{t\ge0}$ satisfies \eqref{eq:cMild}~for $t\in[0,T]$ if and only if 
\begin{equation}
\label{eq:MildV}
V_t\,=\,\Lambda_tX_0 +\int_0^t\Lambda_{t-s}N_\delta(V_s+Y_s)\,ds \qquad\text{ for $t\in[0,T]$}. 
\end{equation}
Let $C_T$ be the space of continuous paths from $[0,T]$ to $D(N)$, and define $U:C_T\rightarrow C_T$ by 
\[
U(v)(t)\,=\,\Lambda_tX_0 + \int_0^t\Lambda_{t-s}N_\delta(v_s+Y_s)\,ds,\qquad v\in C_T,\,\,t\in[0,T].  
\] 
Due to the local boundedness of $(\Lambda_t)_{t\ge0}$ and $N$ on $D(N)$, the local inversion theorem \cite[Lemma 9.2]{DPZ14}~implies that there exists small $T>0$ such that \eqref{eq:MildV}~admits a unique solution $(V_t)_{t\in[0,T]}$, and hence there exists a unique local mild solution $(X'_t)_{t\in[0,T]}$ to \eqref{eq:cSPDE}. 

To extend this solution globally in time, we exploit the fact that $N_\delta$ is Lipschitz to obtain 
\[
e^{-\omega t}\norm*{V_t}_{D(N)}\,\le\,\norm*{X_0} _{D(N)} + \int_0^t \kappa e^{-\omega s}\norm*{Y_s}_{D(N)}\,ds + \int_0^t\kappa e^{-\omega s}\norm*{V_s}_{D(N)}\,ds. 
\]
We may then apply Gr{\"o}nwall's inequality to conclude that $V_t$ cannot blowup in finite time, and we may extend the solution to \eqref{eq:MildV}~to $t\in[0,\infty)$ for $X_0\in D(N)$. 
Extending this solution theory to allow for $X_0\in H$ can be achieved as in \cite[Theorem 7.15]{DPZ14}. 
The same arguments apply to \eqref{eq:SPDE}~for $t<\tau$. 

Similar arguments can be used to prove the existence and uniqueness of solutions in $D(N)$ to \eqref{eq:OU}~with $x,y\in H$. 
Following \cite[Theorem 3.6]{DP04}, we now prove that $x\mapsto X_t'$ is $C^1$ with Fr{\'e}chet derivative $\eta_t$. 
For fixed $x,y\in H$, $t>0$, define for $s\in(0,t)$ the quantity
\[
\rho_s(y)\,\coloneqq\,X'_s(x+y)-X'_s(x) - \eta_s(x)[y]. 
\]
We'll show that $\norm*{\rho_s(y)}_H\le\,c_t\norm*{y}_H^2$ almost surely. 
Using the evolution equations \eqref{eq:cSPDE}~\&~\eqref{eq:Variational}~and the integral remainder form of Taylor's theorem \cite[Theorem 4.C]{Z12}, we have 
\[
\begin{aligned}
\rho_s(y)\,&=\,\int_0^s e^{L(s-r)}\left(N_\delta(X'_r(x+y))-N_\delta(X'_r(x))\right)\,dr \\
&\qquad -\int_0^se^{L(s-r)}DN_\delta(X'_r(x))\eta_r(x)[y]\,dr\\
&=\,\int_0^s e^{L(s-r)}\left[\int_0^1 DN_\delta(\xi X'_r(x+y) + (1-\xi)X_r(x))\,d\xi\right](\rho_r(y)+\eta_r(x)[y])\,dr \\
&\qquad -\int_0^s e^{L(s-r)}DN_\delta(X_r(x))\eta_r[y]\,dr \\
&=\,\int_0^s e^{L(s-r)}\left[\int_0^1 DN_\delta(\xi X'_r(x+y) + (1-\xi)X_r(x))\,d\xi\right]\rho_r(y)\,dr \\
&\qquad +\int_0^s e^{L(s-r)}\left[\int_0^1 DN_\delta(\xi X'_r(x+y) + (1-\xi)X'_r(x)) - DN_\delta(X_r(x))\,d\xi\right]\eta_r[y]\,dr. 
\end{aligned}
\]
Using Assumption \ref{assn:Supp}, for $t>0$ and $s\in(0,t)$, we see that there exists $k_t>0$ such that $\norm*{\eta_s[y]}_{D(N)}\le k_t\norm*{y}_H$. 
One may then use the fact that $N\in C^2(D(N))$ and Gr{\"o}nwall's inequality to obtain 
\[
\begin{aligned}
\norm*{\rho_s(y)}_H\,&\le\,C\norm*{\rho_s(y)}_{D(N)}\,\le\,C\int_0^se^{\omega(s-r)}\,dr\frac{1}{2}\norm*{N_\delta}_{C^2(D(N))} k_t^2\norm*{y}_H^2 \\
&\qquad\qquad\qquad\qquad\qquad\qquad\times\exp\left(\int_0^s\left(\norm*{N_\delta}_{C^1(D(N))}+\norm*{y}_H\right)\,dr\right) \\
&\le\, c_t\norm*{y}_H^2 
\end{aligned}
\]
for a constant $c_t>0$ increasing in $t>0$. 
As a consequence, we see that $\eta_s[y]$ is in fact the Fr{\'e}chet derivative of $x\mapsto X_t'$ for all $t>0$. 
\end{proof}
\end{proposition}

\begin{remark}
Assumption \ref{assn:SPDE}~fails when $B=I$, $L=\Delta$, and $d\ge2$. 
In this case, to make sense of solutions to \eqref{eq:SPDE}~one must perform a renormalization procedure  \cite{DPD03,H13,H14}, as discussed for stochastic reaction-diffusion equations by \cite{BK16}, and for many other classes of SPDEs elsewhere in the literature. 
It is presently unclear as to how much of the following discussion holds for \eqref{eq:SPDE}~when the equation must be renormalized. 
\end{remark}

The following is the main result of this section, proving Theorem A above. 
Henceforth, let $(X'_t)_{t\ge0}$ denote the mild solution to \eqref{eq:cSPDE}, and let $(P_t)_{t\ge0}$ denote the corresponding Markov semigroup, 
\[
P_t:BM(H)\rightarrow BM(H),\qquad P_tf(x)\,\coloneqq\,\ExpOp_x\left[f(X'_t)\right]. 
\]
Let $(X_t)_{t\in[0,\tau)}$ denote the mild solution of \eqref{eq:SPDE}. 

\begin{theorem}
\label{thm:XtPt}
Under Assumption \ref{assn:Supp}, the solution to \eqref{eq:cSPDE}~satisfies Assumptions \ref{assn:Compact}, \ref{assn:Irreducible}~\&~\ref{assn:tau}. 
In particular, the following properties of $(X'_t)_{t\ge0}$ and $(P_t)_{t\ge0}$ hold. 
\begin{enumerate}[(i)]
\item $(P_t)_{t\ge0}$ is irreducible and strong Feller, in the sense of Assumption \ref{assn:Irreducible}. 
\item $(X'_t)_{t\ge0}$ possesses a unique invariant measure $\mu$. 
\item $(P_t)_{t\ge0}$ extends to a strongly continuous semigroup of compact operators on $L^p(H,\mu)$ for $p>1$. 
\end{enumerate}
As a consequence, $(X'_t)_{t\ge0}$ possesses a unique quasi-stationary measure, a unique quasi-ergodic measure, and a unique $Q$-process on $\Gamma_\delta$, which are also the unique quasi-stationary measure, unique quasi-ergodic measure, and unique $Q$-process of the solution to \eqref{eq:SPDE}. 
\end{theorem}

The proof that \eqref{eq:cSPDE}~satisfies Assumption \ref{assn:tau}~follows from \cite[Theorem 12.19]{DPZ14}. 
We prove the remainder of Theorem \ref{thm:XtPt}~in a series of lemmas, beginning with a proof of the stochastic irreducibility and strong Feller property of $(P_t)_{t\ge0}$.  

To prove the strong Feller property and that each $P_t$ is compact on $L^p(H,\mu)$, we make use of results from \cite{PZ95}~and \cite{DPDG02}, respectively. 
To this end we must introduce the Yosida approximations of \eqref{eq:cSPDE}~and \eqref{eq:Variational}, which are essential to the proof of the integration by parts formula necessary for the arguments \cite{DPDG02}~in the case where $N_\delta$ is not a bounded operator on $H$. 
The Yosida approximations of \eqref{eq:SPDE}~and \eqref{eq:Variational}~are constructed by replacing $N_\delta$ with its Yosida approximation $N_{\delta,\alpha}$. 
The Yosida approximation of a dissipative operator is well-defined, and while $N_\delta$ itself is not dissipative, note that $N_\delta-\kappa$ is, where $\kappa=\lip N_\delta$. 
Hence, following the construction in \cite[Appendix D.3]{DPZ14}, we define for small $\alpha>0$ 
\[
J_\alpha(x)\,\coloneqq\,\left(I-\alpha(N_\delta-\kappa I)\right)^{-1}(x),\qquad x\in H, 
\]
and 
\[
(N_{\delta,\alpha}-\kappa I)(x)\,\coloneqq\,(N_\delta-\kappa I)(J_\alpha(x))\,=\,\frac{1}{\alpha}(J_\alpha(x)-x),\qquad x\in H. 
\]
Then, the Yosida approximation of \eqref{eq:cSPDE}~is 
\begin{equation}
\label{eq:SPDEYosida}
dX^\alpha\,=\,\left(LX^\alpha + N_{\delta,\alpha}(X^\alpha)\right)\,dt + \sigma\,dW,
\end{equation}
while the Yosida approximation of \eqref{eq:Variational}~is 
\begin{equation}
\label{eq:VariationalYosida}
\partial_t\eta^\alpha\,=\,\,L\eta_t + DN_{\delta,\alpha}(X_t)\eta_t. 
\end{equation}

\begin{lemma}
\label{lemma:Yosida}
Under Assumption \ref{assn:Supp}, the following hold. 
\begin{enumerate}[(i)]
\item For each $\alpha>0$ the Yosida approximation $N_{\delta,\alpha}:H\rightarrow H$ is Lipschitz, with a Lipschitz constant independent of $\alpha>0$. 
\item The Yosida approximations \eqref{eq:SPDEYosida}~and \eqref{eq:VariationalYosida}~possess unique mild solutions, which we denote $(X^\alpha_t)_{t\ge0}$ and $(\eta_t^\alpha)_{t\ge0}$. 
\item $X^\alpha$ converges to $X'$ almost surely on $[0,\infty)\times H$, and $\eta^\alpha$ converges almost surely to $\eta$ on $[0,\infty)\times H$. 
\end{enumerate}
\begin{proof}
(i)
Taking arbitrary $x,y\in H$ such that $\norm*{J_\alpha(x)-J_\alpha(y)}_H>0$, we observe that 
\begin{equation}
\label{eq:LipJ}
\begin{aligned}
\norm*{J_\alpha(x)-J_\alpha(y)}_{D(N)}\,&=\,\frac{\norm*{J_\alpha(x)-J_\alpha(y)}_{D(N)}}{\norm*{J_\alpha(x)-J_\alpha(y)}_H}\norm*{J_\alpha(x)-J_\alpha(y)}_H\\
&\le\,C\norm*{J_\alpha(x)-J_\alpha(y)}_H\\
&\le\,C\norm*{x-y}_H, 
\end{aligned}
\end{equation}
where $C>0$ is a constant depending on the embedding of $D(N)$ into $H$, and $J_\alpha:H\rightarrow H$ is $1$-Lipshitz by \cite[Proposition D.11(i)]{DPZ14}. 
If $\norm*{J_\alpha(x)-J_\alpha(y)}_H=0$, we may take sequences $(x_n)_{n\in\N},\,(y_n)_{n\in\N}$, such that $x_n\rightarrow x$ and $y_n\rightarrow y$ in $H$, and such that $\norm*{J_\alpha(x_n)-J_\alpha(y_n)}_H>0$. 
Then, from \eqref{eq:LipJ}~we have 
\[
\begin{aligned}
\norm*{J_\alpha(x)-J_\alpha(y)}_{D(N)}\,&=\,\lim_{n\rightarrow\infty}\norm*{J_\alpha(x_n)-J_\alpha(y_n)}_{D(N)}\\
&\le\,\lim_{n\rightarrow\infty}C\norm*{x_n-y_n}_H\\
&=\,C\norm*{x-y}_H. 
\end{aligned}
\]
It follows that, for arbitrary $x,y\in H$, we have 
\[
\begin{aligned}
\norm*{N_{\delta,\alpha}(x)-N_{\delta,\alpha}(y)}_H\,&\le\,C\norm*{N_{\delta}(J_\alpha(x))-N_{\delta}(J_\alpha(y))}_{D(N)} + C\kappa\norm*{J_\alpha(x)-J_\alpha(y)}_{D(N)}\\
&\le\,2C\kappa\norm*{J_\alpha(x)-J_\alpha(y)}_{D(N)}\\
&\le\,2C^2\kappa\norm*{x-y}_H,
\end{aligned}
\]
concluding the proof of (i). 

(ii)
The Yosida approximations $N_{\delta,\alpha}$ are Lipschitz and bounded on $H$, so that the existence and uniqueness of solutions to \eqref{eq:SPDEYosida}~\&~\eqref{eq:VariationalYosida}~follow from the same arguments as in Proposition \ref{thm:Solution}. 

(iii) 
We sketch the proof of the convergence of $(X_t^\alpha)_{t\ge0}$ to $(X'_t)_{t\ge0}$. 
Observe that 
\begin{equation}
\label{eq:YosidaControl}
\begin{aligned}
\norm*{X'_t-X_t^\alpha}_{D(N)}\,&=\,\norm*{\int_0^t\Lambda_{t-s}\left(N_\delta(X'_s)-N_{\delta,\alpha}(X_s^\alpha)\right)\,ds}_{D(N)}\\
&\le\, \int_0^te^{-\omega(t-s)}\norm*{N_\delta(X'_s)-N_\delta(X_s^\alpha)}_{D(N)}\,ds \\
&\qquad+\int_0^t e^{-\omega(t-s)}\norm*{N_\delta(X_s^\alpha)-N_{\delta,\alpha}(X_s^\alpha)}_{D(N)}\,ds. 
\end{aligned}
\end{equation}
Using the fact that $N_{\delta,\alpha}(x)=N_\delta(J_\alpha(x))$, by the proof of \cite[Proposition D.11(iii)]{DPZ14}~and the Lipschitz property of $N_\delta$, for $x\in D(N)$ we have 
\[
\norm*{N_\delta(x)-N_{\delta,\alpha}(x)}_{D(N)}\,\le\,\kappa\norm*{x-J_\alpha(x)}\,\le\,\kappa\alpha\norm*{(N_\delta-\kappa)x}_{D(N)}. 
\]
Therefore from \eqref{eq:YosidaControl}~we may apply Gr{\"o}nwall's inequality and take $\alpha\rightarrow0$ to conclude the result. 
A similar argument applies for the convergence of $(\eta_t^\alpha)_{t\ge0}$ to $(\eta_t)_{t\ge0}$. 
\end{proof}
\end{lemma}

\begin{lemma}
\label{lemma:Irreducible}
Assumption \ref{assn:Supp}~implies that $(P_t)_{t\ge0}$ is irreducible in the sense of \eqref{eq:Irreducible}. 
Moreover, $(P_t)_{t\ge0}$ is a strong Feller semigroup, in the sense that for any bounded measurable function $f:H\rightarrow\R$, the function $P_tf:H\rightarrow\R$ is bounded and continuous. 
\begin{proof}
Stochastic irreducibility follows from Propositions 2.8 \&~2.11 of \cite{M93}. 
To prove that $(P_t)_{t\ge0}$, we apply a result of \cite{PZ95}~to $(P_t^\alpha)_{t\ge0}$, and then take $\alpha\rightarrow0$. 
One can verify that assumptions (A.1), (A.2), \&~(A.3) of \cite{PZ95}~are satisfied for the Yosida approximation \eqref{eq:SPDEYosida}. 
In particular, by Lemma \ref{lemma:Yosida}~above, $N_{\delta,\alpha}:H\rightarrow H$ is Lipschitz with Lipschitz constant $2C^2\kappa$. 
Applying \cite[Theorem 1.2]{PZ95}, it follows that for arbitrary $\alpha>0$, $t>0$, $x,y\in H$, and bounded measurable $f:H\rightarrow\R$, we have 
\begin{equation}
\label{eq:YosidaFeller}
\abs*{P_t^\alpha f(x)- P_t^\alpha f(y)}\,\le\,c_t\norm*{f}_\infty\norm*{x-y}, 
\end{equation}
where in our setting $c_t=(t\wedge t_0)^{-1/2}e^{-\omega t_0}$ for a time $t_0>0$ such that 
\[
t_0^2e^{-2\omega t_0}+\int_0^{t_0}\norm*{\Lambda_t}_{HS}^2\,dt\,\le\,\frac{2}{9}(2C^2\kappa)^{-2}.
\] 
Since $c_t$ is independent of $\alpha>0$, we may take $\alpha\rightarrow0$ in \eqref{eq:YosidaFeller}~to conclude that 
\[
\abs*{P_t f(x)- P_t f(y)}\,\le\,c_t\norm*{f}_\infty\norm*{x-y}. 
\]
Hence, for any bounded measurable $f:H\rightarrow\R$, we have that $P_tf:H\rightarrow\R$ is continuous. 
\end{proof}
\end{lemma}

\begin{lemma}
\label{thm:CutoffInvariant}
Assumption \ref{assn:Supp}~implies that \eqref{eq:cSPDE}~possesses a unique invariant measure $\mu$ on $H$. 
Moreover, $\mu$ satisfies a logarithmic Sobolev inequality. 
\begin{proof}
To prove the existence and uniqueness of an invariant measure of \eqref{eq:cSPDE}, we apply \cite[Theorem 2.5]{GM05}. 
To do so, we must verify that \eqref{eq:cSPDE}~satisfies  
\begin{itemize}
\item \cite[Hypothesis 2.1]{GM05}~$(P_t)_{t\ge0}$ is irreducible and strong Feller. 
\item \cite[Hypothesis 2.2]{GM05}~For each $r>0$, there exists some $t_0>0$ and a compact set $K\subset D(N)$ such that $\inf\left\{P_{t_0}1_K(x)\,:\,\norm*{x}_{D(N)}\le r\right\}>0$. 
\end{itemize}

Lemma \ref{lemma:Irreducible}~above implies that \cite[Hypothesis 2.1]{GM05}~is satisfied. 
To prove that $(P_t)_{t\ge0}$ satisfies \cite[Hypothesis 2.2]{GM05}, we borrow ideas from \cite{GM05,MS99}, filling in the details of their arguments for the case where $N$ is only defined on $D(N)\subsetneq H$. 
We first show that for arbitrary $t>0$ and $f\in C_b(H)$, $P_tf$ is weakly sequentially continuous. 
Hence take $(x_n)_{n\in\N}\subset H$ converging weakly to $x\in H$, and observe that 
\[
\begin{aligned}
\ExpOp\left[\norm*{X_t(x_n)-X_t(x)}_{D(N)}\right]\,&\le\,\norm*{\Lambda_t(x_n-x)}_{D(N)} \\
&\qquad\qquad+ \ExpOp\left[\norm*{\int_0^t\Lambda_{t-s}(N_\delta(X_s(x_n))-N_\delta(X_s(x)))\,ds}_{D(N)}\right] \\
&\le\, \norm*{\Lambda_t(x_n-x)}_{D(N)} \\
&\qquad\qquad+ \ExpOp\left[\int_0^t\norm*{\Lambda_{(t-s)}}_{D(N)}\norm*{N_\delta(X_s(x_n))-N_\delta(X_s(x))}_{D(N)}\,ds\right]\\
&\le\, \norm*{\Lambda_t(x_n-x)}_{D(N)} + \ExpOp\left[\int_0^te^{-\omega(t-s)}\kappa\norm*{X_s(x_n)-X_s(x)}_{D(N)}\,ds\right]. 
\end{aligned}
\]
Define $a_t^n\coloneqq\norm*{\Lambda_t(x_n-x)}_{D(N)}$. 
By Assumption \ref{assn:PDE}, $\Lambda_t$ is compact, and therefore completely continuous. 
Hence for each $t>0$ and any $\epsilon>0$, there exists $N\in\N$ such that if $n\ge N$ then $a_t^n\le\epsilon$. 
Applying Tonelli's theorem, we therefore find that for large enough $n\in\N$ we have 
\[
e^{\omega t}\ExpOp\left[\norm*{X_t(x_n)-X_t(x)}_{D(N)}\right]\,\le\,\epsilon + \int_0^t\kappa e^{\omega s}\ExpOp\left[\norm*{X_s(x_n)-X_s(x)}_{D(N)}\right]\,ds. 
\]
Applying Gr{\"o}nwall's inequality, we find that for any $t>0$, 
\[
c\ExpOp\left[\norm*{X_t(x_n)-X_t(x)}_H\right]\,\le\,\ExpOp\left[\norm*{X_t(x_n)-X_t(x)}_{D(N)}\right]\,\le\,\epsilon e^{-\omega t} + \int_0^t\epsilon \kappa^s\,ds, 
\]
which converges to zero as $n\rightarrow\infty$, where $c>0$ is a constant arising from the embedding of $D(N)$ in $H$. 
From this, one may readily see that for any bounded continuous $f:H\rightarrow\R$ we have that $f(X_t(x_n))\rightarrow f(X_t(x))$ in probability as $n\rightarrow\infty$. 
By the dominated convergence theorem, $P_tf(x_n)$ must also converge to $P_tf(x)$ as $n\rightarrow\infty$, completing the proof of the claim that $P_tf$ is weakly sequentially continuous for $f\in C_b(H)$. 

Now, fix $r>0$, $t_0>0$, and $\epsilon>0$, and define $B_r\coloneqq\left\{x\in D(N)\,:\,\norm*{x}_{D(N)}\le r\right\}$. 
Let $\overline{B}_r$ be the closure of $B_r$ in the weak topology of $H$. 
Viewing $P_t1_\cdot(x)=\ExpOp_x[X_t\in\cdot]$ as a measure on $H$, by the definition of the narrow topology we have that $x\mapsto P_t1_\cdot(x)$ is a continuous map from $\overline{B}_r$ equipped with the weak topology to the space of probability measures on $H$ equipped with the narrow topology. 
Hence, by the compactness of $\overline{B}_r$ in $H$, we see that $\mathcal{M}\coloneqq\left\{P_t1_\cdot(x)\,:\,x\in\overline{B}_r\right\}$ is narrowly compact. 
As a consequence, by the Prokhorov theorem \cite[Theorem 8.6.2]{B00}~the family of measures $\mathcal{M}$ is tight\footnote{Here, we say that a sequence of measures $\mu_n$ converges ``in the narrow topology'' if $\mu_n(f)\rightarrow\mu(f)$ for all bounded continuous $f$. 
This is equivalent to what \cite{B00}~refers to as weak convergence.}, which is equivalent to saying that \cite[Hypothesis 2.2]{GM05}~holds. 

Finally, we note that there exist constants $c_1,c_2>0$ such that 
\[
\begin{aligned}
\ExpOp_x\left[\norm*{X'_t}_{D(N)}\right]\,&\le\,e^{-\omega t}\ExpOp_x\left[\norm*{X_0}_{D(N)}\right] + \int_0^t e^{-\omega(t-s)}\ExpOp_x\left[\norm*{N_\delta(X'_s)}_{D(N)}\right]\,ds \\
&\qquad\qquad\qquad\qquad\qquad\qquad\qquad + \sigma\ExpOp_x\left[\norm*{Y_t}_{D(N)}\right]\\
&\le\,c_1e^{-\omega t} + c_2(1-e^{-\omega t}), 
\end{aligned}
\]
and hence we may apply \cite[Theorem 2.5]{GM05}~to conclude the existence and uniqueness of an invariant measure $\mu$ for \eqref{eq:cSPDE}. 

We now sketch the proof of the logarithmic Sobolev inequality, following ideas from \cite[Proposition 3.30]{DP04}. 
First, note that Assumption \ref{assn:Supp}~implies  
\[
\norm*{DX'_t}_{D(N)}\,\le\,e^{-(\omega-\kappa)t}, 
\]
with $\omega>\kappa$, and therefore for any $f\in C^1_b(H)$ we have for some constant $C_f>0$ that 
\[
\begin{aligned}
\norm*{DP_t(f^2)}^2\,&=\,\norm*{2\ExpOp_x\left[f(X'_t)Df(X'_t)DX'_t\right]}^2\\
&\le\,2C_fe^{-(\omega-\kappa)t}. 
\end{aligned}
\]
Consequently, taking into account $\mu$-invariance and the integration by parts formula \cite[Theorem 3.6]{DPDG02}, 
\[
\begin{aligned}
\partial_t\int P_tf^2\ln P_tf^2\,d\mu\,&=\,\int\mathcal{L} P_tf^2\ln P_tf^2\,d\mu + \int\mathcal{L}P_tf^2\,d\mu\\
&=\,-\frac{1}{2}\int\frac{1}{P_tf^2}\norm*{BDP_tf^2}^2\,d\mu \\
&\ge\,-C_fe^{-(\omega-\kappa)t}\int P_t\norm*{BDf}^2\,d\mu\,=\,-C_f\mu\left(\norm*{BDf}^2\right)e^{-(\omega-\kappa)t}. 
\end{aligned}
\]
Integrating the above over a finite time interval $[0,t]$, we obtain 
\[
\int P_sf^2\ln P_sf^2 \,d\mu\big\rvert_0^t\,\ge\,(\omega-\kappa)^{-1}(1-e^{-(\omega-\kappa)t})\mu\left(\norm*{BDf}^2\right), 
\]
and taking $t\rightarrow\infty$ yields the logarithmic Sobolev inequality  
\[
\mu(f^2)\ln\mu(f^2) - \mu(f^2\ln f^2)\,\ge\,(\omega-\kappa)^{-1}\mu\left(\norm*{BDf}^2\right). 
\] 
\end{proof}
\end{lemma}

Finally, to apply the results of \cite{DPDG02}, we must prove a result relating the invariant measure $\mu$ of \eqref{eq:cSPDE}~to the invariant measure of \eqref{eq:OU}. 

\begin{lemma}
\label{lemma:OU}
Under Assumption \ref{assn:Supp}, the Ornstein-Uhlenbeck process governed by \eqref{eq:OU}~possesses a unique invariant measure $\nu$, which is a Gaussian measure with mean zero and covariance operator 
\[
Q_\infty\,\coloneqq\,\frac{\sigma^2}{2}\int_0^\infty\Lambda_tBB^*\Lambda_t^*\,dt. 
\]
Moreover, the density $\rho\coloneqq d\mu/d\nu$ exists, and $\ln\rho\,\in\,W^{1,2}(H,\mu)$. 
\begin{proof}
The existence of the invariant measure $\nu$ follows from \cite[Theorem 2.34]{DP04}. 
To prove the existence of $\rho$ and that $\ln\rho\in W^{1,2}(H,\mu)$, first observe that 
\begin{equation}
\label{eq:Dlnrho} 
\int_H \norm*{N_\delta(x)}_{D(N)}^2\,\mu(dx)\,\le\,\sup_{x\in\Gamma_\delta}\norm*{N(x)}_{D(N)}^2\mu(\Gamma_\delta)\,<\,\infty.  
\end{equation}
Letting $\mu_\alpha$ be the unique invariant measure of the Yosida approximation \eqref{eq:SPDEYosida}, the existence of $\rho_\alpha=d\mu_\alpha/d\nu$ follows from \cite{BDPR96}. 
Moreover, \cite[Claim 3]{BDPR96}~implies that $D\ln\rho_\alpha= N_{\delta,\alpha}$, where $D\ln\rho_\alpha$ denotes the Fr{\'e}chet derivative of $\ln\rho_\alpha$, and hence  
\[
\int_H \norm*{D\ln\rho_\alpha(x)}^2\,\mu(dx)\,=\,\int_H\norm*{N_{\delta,\alpha}(x)}^2\,\mu(dx)\,<\,\infty, 
\]
implying that $\ln\rho_\alpha\in W^{1,2}(H,\mu)$. 
For $\alpha>0$, let $(P_t^\alpha)_{t\ge0}$ denote the Markov semigrop of $(X_t^\alpha)_{t\ge0}$. 
Since $X_t^\alpha$ converges to $X_t$ in $D(N)$ for each $t>0$ as $\alpha\rightarrow0$, by the dominated convergence theorem we have 
\[
\lim_{\alpha\rightarrow0}\mu_\alpha(A)\,=\,\lim_{\alpha\rightarrow0}\lim_{t\rightarrow\infty}P_t^\alpha1_A(x)\,=\,\lim_{t\rightarrow\infty}P_t1_A(x)\,=\,\mu(A). 
\]
Hence by the Vitali-Hahn-Sack Theorem \cite[Chapter II.2]{Y96}, $\rho=d\mu/d\nu$ exists. 
By \eqref{eq:Dlnrho}~we may apply the dominated convergence theorem to conclude that $\ln\rho\in W^{1,2}(H,\mu)$. 
\end{proof}
\end{lemma}

We are now able to prove the fourth statement in Theorem \ref{thm:XtPt}

\begin{lemma}
\label{prop:Compact}
Let Assumption \ref{assn:Supp}~hold. 
For each $p\in[1,\infty)$, $(P_t)_{t\ge0}$ extends to a strongly continuous semigroup of compact operators on $L^p(H,\mu)$. 
\begin{proof}
Follows from \cite[Theorem 5.1]{DPDG02}. 
\end{proof}
\end{lemma}

We have therefore proven Theorem \ref{thm:XtPt}~for the solution to \eqref{eq:SPDE}~with $N$ replaced by $N_\delta$, and hence this system with a cutoff nonlinearity possesses a unique quasi-stationary measure, a unique quasi-ergodic measure, and a unique $Q$-process in $\Gamma_\delta$ defined $\mu$-almost surely in $\Gamma_\delta$. 
However, our interest is ultimately in \eqref{eq:SPDE}~with a non-cutoff nonlinearity. 
The following result demonstrates that the quasi-stationary distribution, quasi-ergodic distribution, and $Q$-process of \eqref{eq:SPDE}~with cutoff nonlinearity in $\Gamma_\delta$ are also the quasi-stationary measure, quasi-ergodic measure, and $Q$-process of \eqref{eq:SPDE}~with non-cutoff nonlinearity. 

\begin{proposition}
\label{prop:CutoffEquivalence}
Let Assumption \ref{assn:Supp}~hold, and let $\alpha$ and $\beta$ be the unique quasi-stationary and quasi-ergodic measure of \eqref{eq:cSPDE}, guaranteed to exist by Theorems \ref{thm:XtPt}~\&~\ref{thm:quasiErgodic}. 
Then, $\alpha$ and $\beta$ are the unique quasi-stationary and quasi-ergodic measure of \eqref{eq:SPDE}. 
Moreover, \eqref{eq:SPDE}~admits a unique $Q$-process in $\Gamma_\delta$, equal to the $Q$-process of \eqref{eq:cSPDE}. 
\begin{proof}
For $t<\tau$ we of course have $N(X_t)=N_\delta(X_t)$ and $N(X'_t)=N_\delta(X'_t)$. 
Since $(X_t)_{t\ge0}$ and $(X'_t)_{t\ge0}$ are driven by the same Wiener process, it follows that
\[
\norm*{X_t-X_t'}\,\le\,\int_0^t e^{-\omega(t-s)}\kappa\norm*{X_s-X'_s}\,ds, 
\]
where $\kappa>0$ is the Lipschitz constant of $N_\delta$. 
Gr{\"o}nwall's inequality then implies $\norm*{X_t-X'_t}=0$ for $t<\tau$. 
In particular, $\pr_x[X_t\in\cdot\,|\,t<\tau]=\pr_x[X'_t\in\cdot\,|\,t<\tau]$. 
\end{proof}
\end{proposition}

\section{Applications}
\label{sec:Applications}
In this section, we sketch how the results of Section \ref{sec:Examples}~may be used to study the effects of noise on the behaviour of metastable patterns in SPDEs. 
We first study how noise may affect the average ``position'' of a pattern. 
Here, the position of a pattern is rigorously defined using a generalization of the \emph{isochron map}, introduced by Winfree \cite{W74}~to define the phase of stochastic oscillators in biological systems. 
Loosely speaking, the isochron map uniquely projects a given point in the vicinity of a deterministic stable invariant manifold $\Gamma$ to a point on $\Gamma$, in a manner that is consistent with the deterministic dynamics. 
The generalization of the isochron map used here allows us to define the \emph{isochronal phase}~of patterns such as stochastic travelling waves and stochastic spiral waves, giving a precise notion of the position of such a stochastically perturbed pattern. 
This is done in Section \ref{sec:Isochron}. 
Using the isochron map, in Section \ref{sec:Frequency}~we sketch how one may obtain results on noise-induced deviations in the \emph{velocity}~of a spatiotemporal pattern. 
To do so we restrict our attention to the case where the pattern $\Gamma$ is a limit cycle. 

\subsection{Quasi-Ergodicity of the Isochronal Phase} 
\label{sec:Isochron}
Let $(X_t)_{t\ge0}$ denote the unique mild solution of the following stochastic reaction-diffusion system 
\begin{equation}
\label{eq:ReacDiff}
dX\,=\,\left((\Delta-a)X + N(X)\right)\,dt + \sigma\,dW, 
\end{equation}
where $W$ is a cylindrical Wiener process, $a,\sigma>0$ are parameters, $N:D(N)\subset H\rightarrow H$ is a vector of polynomials, and $H=L^2(O,\R^m)$ is the solution space of \eqref{eq:ReacDiff}~with a bounded spatial domain $O$. 

\begin{assumption}
\label{assn:ReacDiff}
Let $(t,x)\mapsto\phi_t(x)$ denote the solution map of \eqref{eq:ReacDiff}~with $\sigma=0$. 
There exists a stable normally hyperbolic $C^1$ invariant manifold $\Gamma$ of \eqref{eq:ReacDiff}~which is compact, smooth, and finite dimensional. 
Moreover, 
\begin{enumerate}[(a)]
\item $\phi_t$ is bi-Lipschitz in $\Gamma_\delta$, the $\delta$-neighbourhood of $\Gamma$ defined in the topology of $D(N)$. 
In particular, the dynamics of $\phi_t$ are not chaotic near $\Gamma$, 
\item $D\phi_t(\gamma)$ is invertible with an inverse bounded uniformly in $\gamma\in\Gamma$, 
\item $\Gamma$ can be parameterized by some $\mathcal{S}\subset\R^{m'}$ as $\Gamma=\{\gamma_\alpha\}_{\alpha\in\mathcal{S}}$, and $\alpha\mapsto\gamma_\alpha$ is $C^1$ with $D\gamma_\alpha$ being invertible with an inverse bounded uniformly in $\alpha$. 
\end{enumerate}
\end{assumption}

To study the long-term behaviour of $(X_t)_{t\ge0}$ in the vicinity of $\Gamma_\delta$, we reduce this process to a finite dimensional process on $\Gamma$ in a manner that is consistent with the deterministic dynamics. 

\begin{lemma}
Let Assumption \ref{assn:ReacDiff}~hold. 
Then, for each $x\in\Gamma_\delta$ there exists a unique point $\pi(x)\in\Gamma$ such that 
\begin{equation}
\label{eq:Isochron}
\norm*{\phi_t(x) - \phi_t(\pi(x))}\,\xrightarrow[t\rightarrow\infty]{}\,0. 
\end{equation} 
The map $\pi:\Gamma_\delta\rightarrow\Gamma$ is referred to as the \emph{isochron map}~of $\Gamma$. 
Moreover, if $\phi_t$ is $C^2$, then $\pi$ is $C^2$ (in the topology of $D(N)$). 
\begin{proof}
A simple proof follows from Cantor's Intersection Theorem, as in \cite{A21}, or, under slightly stronger assumptions, from a fixed point argument as in \cite[Theorem 3.1]{AM22}. 
\end{proof}
\end{lemma}

Remarkably, even though \eqref{eq:SPDE}~only admits a mild solution when the noise is not trace class, the isochron map possesses certain regularizing properties that allow one to prove a strong It{\^o}~formula. 

\begin{lemma}
Let Assumption \ref{assn:ReacDiff}~hold. 
Let $V=\Delta-a+N$, and let $D\pi$ and $D^2\pi$ denote the first and second Fr{\'e}chet derivatives of $\pi$ in the topology of $D(N)$, respectively. 
Then, for some orthonormal basis $\{e_k\}_{k\in\N}$ of $H$ it holds that  
\begin{equation}
\label{eq:ItoLemma}
\pi(X_{t\wedge\tau})-\pi(X_0)\,=\, \int_0^{t\wedge\tau} D\pi(X_s)V(X_s) + \frac{\sigma^2}{2}\sum_{k\in\N}D^2\pi(X_s)[Be_k,Be_k] \,ds + \int_0^{t\wedge\tau} D\pi(X_s)B\,dW_s, 
\end{equation}
and all of the above integrals are well defined. 
\begin{proof}
See \cite[Theorem 3.5]{AM22}. 
\end{proof}
\end{lemma}

Now, so long as \eqref{eq:SPDE}~satisfies Assumption \ref{assn:Supp}, we may apply Theorem A to conclude that \eqref{eq:SPDE}~admits a unique quasi-ergodic measure $\beta$ in $\Gamma_\delta$. 
We may then define a measure on $\Gamma$ via pushforward under the isochron map, 
\[
\pi_*\beta(A)\,\coloneqq\,\beta\left(\pi^{-1}(A)\right). 
\]
Using the change of variables formula and Theorems \ref{thm:quasiErgodic}~\&~\ref{thm:XtPt}, we have the following result, providing a measure that indicates how the noise affects the average position of our stochastically perturbed pattern. 
This may be compared with \cite[Theorem C]{AM22}, obtained using different methods.  

\begin{theorem}
\label{thm:A}
Let Assumption \ref{assn:ReacDiff}~hold. 
For any $p\ge1$, $\epsilon>0$, and bounded $g\in L^2(\Gamma,\pi_*\beta)$, letting $\pi_t\coloneqq\pi(X_t)$ for $t<\tau$ we have 
\begin{equation}
\lim_{t\rightarrow\infty}\pr_x\left[\abs*{\frac{1}{t}\int_0^tg(\pi_s)\,ds - (\pi_*\beta)(g)}^p>\varepsilon\,\big|\,t<\tau\right]\,=\,0. 
\end{equation}
\end{theorem}

\subsection{Quasi-Asymptotic Frequencies of Stochastic Oscillators}
\label{sec:Frequency}
Now, assume that $\Gamma$ is a periodic (in time) solution of \eqref{eq:PDE}. 
In this case, we refer to a $\Gamma$-like solutions as a \emph{stochastic oscillator}. 
For instance, $\Gamma$ could be a family of periodic solutions corresponding to a limit cycle solution of \eqref{eq:PDE}~without diffusion (\emph{i.e.}~of the ODE $\dot{x}=N(x)$), or a stable pulse-type travelling wave solution of \eqref{eq:PDE}~when $O$ is a periodic spatial domain, as in \cite{AK15}. 
Let the period of $\Gamma$ under the flow of \eqref{eq:PDE}~be $T>0$. 
Specifically, for $(t,x)\in[0,\infty)\times H$ we have 
\[
\begin{aligned}
\gamma_{t+T}\,=\,\gamma_t,\quad\text{ and }\quad\phi_s(\gamma_t)\,=\,\gamma_{s+t},\qquad\forall\,t\in\R,\,\,s\ge0.  
\end{aligned}
\]

For a sufficiently small neighbourhood $\Gamma_\delta$ of $\Gamma$, define the isochron map $\pi:\Gamma_\delta\rightarrow[0,T]$ as in \eqref{eq:Isochron}. 
The stochastic integral in \eqref{eq:ItoLemma}, being a finite dimensional local Markov process with linear quadratic variation, has a time average equal to zero -- see for instance \cite{vZ00}~for details. 
Applying the second statement of Theorem \ref{thm:quasiErgodic}, we find that 
\begin{equation}
\label{eq:quasiFrequency} 
c_\sigma\,\coloneqq\,\lim_{t\rightarrow\infty}\ExpOp_x\left[\frac{1}{t}\pi(X_t)\,|\,t<\tau\right]\,=\,\int_E \pi'(x)V(x) + \frac{\sigma^2}{2}\sum_{k\in\N} \pi''(x)[Be_k,Be_k]\,\beta(dx) 
\end{equation}
exists as a deterministic limit. 
We refer to this $c_\sigma$ as the \emph{quasi-asymptotic frequency}~of \eqref{eq:SPDE}~in $\Gamma_\delta$. 

Based on the properties of the isochron map, it can be shown as in \cite{A22}~that 
\[
c_\sigma\,=\,c_0 + \sigma^2(b_0+b_\sigma), 
\]
where $c_0=1/T$ is the frequency of $\Gamma$ under the deterministic flow of \eqref{eq:PDE}, and $b_0$ is a constant independent of the noise amplitude. 
Moreover, 
\[
b_\sigma\,=\, \frac{1}{2}\int_{\Gamma_\delta} \sum_{k\in\N}\pi''(x)\left[Be_k,Be_k\right]\,\delta_\sigma(dx), 
\]
where $\delta_\sigma\coloneqq\beta-\eta$ is the $\sigma$-dependent difference measure of the quasi-ergodic measure $\beta$ and the deterministic invariant measure $\eta$ of $\Gamma$ under the deterministic flow of \eqref{eq:PDE}. 
Thus, we see that the quasi-asymptotic frequency $c_\sigma$ depends ``almost'' quadratically on $\sigma$ in some interval $[\sigma_0,\sigma_1]\subset[0,\infty)$ if and only if $db_\sigma/d\sigma$ is ``almost'' zero for all $\sigma\in[\sigma_0,\sigma_1]$. 
This might be considered a refinement of the results of \cite{GPS18}, who suggest that this asymptotic frequency should depend quadratically on $\sigma>0$. 
The approach taken here allows for a straightforward of their results to the infinite dimensional setting.

\bibliographystyle{amsplain}
\bibliography{EJP_ZPAdams_QuasiErgodicity_of_SPDE}

\providecommand{\bysame}{\leavevmode\hbox to3em{\hrulefill}\thinspace}
\providecommand{\MR}{\relax\ifhmode\unskip\space\fi MR }
\providecommand{\MRhref}[2]{%
  \href{http://www.ams.org/mathscinet-getitem?mr=#1}{#2}
}
\providecommand{\href}[2]{#2}
\begin{thebibliography}{10}

\bibitem{A22}
Zachary~P Adams, \emph{The asymptotic frequency of stochastic oscillators},
  SIAM Journal on Applied Dynamical Systems \textbf{22} (2023), no.~1,
  311--338.

\bibitem{A21}
Zachary~P. Adams, \emph{Existence, regularity, and a strong it{\^o}~formula for
  the isochronal phase of spde}, arXiv preprint arXiv:2109.04515v3 (2023).

\bibitem{AM22}
Zachary~P Adams and James MacLaurin, \emph{The isochronal phase of stochastic
  pde and integral equations: Metastability and other properties}, arXiv
  preprint arXiv:2210.10681 (2022).

\bibitem{AB80}
Charalambos~D Aliprantis and Owen Burkinshaw, \emph{Positive compact operators
  on banach lattices}, Mathematische Zeitschrift \textbf{174} (1980), 289--298.

\bibitem{AK15}
Gianni Arioli and Hans Koch, \emph{Existence and stability of traveling pulse
  solutions of the fitzhugh--nagumo equation}, Nonlinear Analysis: Theory,
  Methods \& Applications \textbf{113} (2015), 51--70.

\bibitem{BE93}
Markus B{\"a}r and Markus Eiswirth, \emph{Turbulence due to spiral breakup in a
  continuous excitable medium}, Physical Review E \textbf{48} (1993), no.~3,
  R1635.

\bibitem{B15}
Florent Barret, \emph{Sharp asymptotics of metastable transition times for one
  dimensional spdes}, Annales de l'IHP Probabilit{\'e}s et statistiques,
  vol.~51, 2015, pp.~129--166.

\bibitem{BLZ98}
Peter~W Bates, Kening Lu, and Chongchun Zeng, \emph{Existence and persistence
  of invariant manifolds for semiflows in banach space}, vol. 645, American
  Mathematical Soc., 1998.

\bibitem{BFR17}
Andr{\'a}s B{\'a}tkai, M~Kramar Fijav{\v{z}}, and Abdelaziz Rhandi,
  \emph{Positive operator semigroups}, Operator Theory: advances and
  applications \textbf{257} (2017).

\bibitem{BG13}
Nils Berglund and Barbara Gentz, \emph{Sharp estimates for metastable lifetimes
  in parabolic spdes: Kramers' law and beyond}, Electronic Journal of
  Probability \textbf{18} (2013), 1--58.

\bibitem{BK16}
Nils Berglund and Christian Kuehn, \emph{Regularity structures and
  renormalisation of fitzhugh--nagumo spdes in three space dimensions},
  Electronic Journal of Probability \textbf{21} (2016), 1--48.

\bibitem{BL08}
Wolf-J{\"u}rgen Beyn and Jens Lorenz, \emph{Nonlinear stability of rotating
  patterns}, Dynamics of Partial Differential Equations \textbf{5} (2008),
  no.~4, 349--400.

\bibitem{BG16}
Alessandra Bianchi and Alexandre Gaudilliere, \emph{Metastable states,
  quasi-stationary distributions and soft measures}, Stochastic Processes and
  their Applications \textbf{126} (2016), no.~6, 1622--1680.

\bibitem{BGM20}
Alessandra Bianchi, Alexandre Gaudilli{\`e}re, and Paolo Milanesi, \emph{On
  soft capacities, quasi-stationary distributions and the pathwise approach to
  metastability}, Journal of Statistical Physics \textbf{181} (2020), no.~3,
  1052--1086.

\bibitem{B57}
Garrett Birkhoff, \emph{Extensions of jentzsch's theorem}, Transactions of the
  American Mathematical Society \textbf{85} (1957), no.~1, 219--227.

\bibitem{BDPR96}
Vladimir~I Bogachev, Giuseppe Da~Prato, and Michael R{\"o}ckner,
  \emph{Regularity of invariant measures for a class of perturbed
  ornstein-uhlenbeck operators}, Nonlinear Differential Equations and
  Applications NoDEA \textbf{3} (1996), no.~2, 261--268.

\bibitem{B00}
Vladimir~Igorevich Bogachev, \emph{Measure theory}, vol.~2, Springer, 2007.

\bibitem{BW12}
Paul~C Bressloff and Matthew~A Webber, \emph{Front propagation in stochastic
  neural fields}, SIAM Journal on Applied Dynamical Systems \textbf{11} (2012),
  no.~2, 708--740.

\bibitem{BG09}
Thomas Butler and Nigel Goldenfeld, \emph{Robust ecological pattern formation
  induced by demographic noise}, Physical Review E \textbf{80} (2009), no.~3,
  030902.

\bibitem{CLMR21}
Matheus~M Castro, Jeroen~SW Lamb, Guillermo~Olic{\'o}n M{\'e}ndez, and Martin
  Rasmussen, \emph{Existence and uniqueness of quasi-stationary and
  quasi-ergodic measures for absorbing markov processes: a banach lattice
  approach}, arXiv preprint arXiv:2111.13791 (2021).

\bibitem{CV16}
Nicolas Champagnat and Denis Villemonais, \emph{Exponential convergence to
  quasi-stationary distribution and $q$-process}, Probability Theory and
  Related Fields \textbf{164} (2016), no.~1, 243--283.

\bibitem{CV17a}
\bysame, \emph{General criteria for the study of quasi-stationarity}, arXiv
  preprint arXiv:1712.08092 (2017).

\bibitem{CV17b}
\bysame, \emph{Uniform convergence to the $ q $-process}, Electronic
  Communications in Probability \textbf{22} (2017), 1--7.

\bibitem{CV21}
\bysame, \emph{Lyapunov criteria for uniform convergence of conditional
  distributions of absorbed markov processes}, Stochastic Processes and their
  Applications \textbf{135} (2021), 51--74.

\bibitem{C85}
Kai~Lai Chung, \emph{Doubly-feller process with multiplicative functional},
  Seminar on stochastic processes, 1985, Springer, 1985, pp.~63--78.

\bibitem{CMSM13}
Pierre Collet, Servet Mart{\'\i}nez, and Jaime San~Mart{\'\i}n,
  \emph{Quasi-stationary distributions: Markov chains, diffusions and dynamical
  systems}, vol.~1, Springer, 2013.

\bibitem{CRM16}
George~WA Constable, Tim Rogers, Alan~J McKane, and Corina~E Tarnita,
  \emph{Demographic noise can reverse the direction of deterministic
  selection}, Proceedings of the National Academy of Sciences \textbf{113}
  (2016), no.~32, E4745--E4754.

\bibitem{C05}
Stephen Coombes, \emph{Waves, bumps, and patterns in neural field theories},
  Biological cybernetics \textbf{93} (2005), no.~2, 91--108.

\bibitem{DP04}
Giuseppe Da~Prato, \emph{Kolmogorov equations for stochastic pdes}, Springer
  Science \& Business Media, 2004.

\bibitem{DPD03}
Giuseppe Da~Prato and Arnaud Debussche, \emph{Strong solutions to the
  stochastic quantization equations}, The Annals of Probability \textbf{31}
  (2003), no.~4, 1900--1916.

\bibitem{DPDG02}
Giuseppe Da~Prato, Arnaud Debussche, and Beniamin Goldys, \emph{Some properties
  of invariant measures of non symmetric dissipative stochastic systems},
  Probability theory and related fields \textbf{123} (2002), 355--380.

\bibitem{DPZ96}
Giuseppe Da~Prato and Jerzy Zabczyk, \emph{Ergodicity for infinite dimensional
  systems}, vol. 229, Cambridge University Press, 1996.

\bibitem{DPZ14}
\bysame, \emph{Stochastic equations in infinite dimensions}, Cambridge
  university press, 2014.

\bibitem{D84}
Edward~Brian Davies and Barry Simon, \emph{Ultracontractivity and the heat
  kernel for schr{\"o}dinger operators and dirichlet laplacians}, Journal of
  Functional Analysis \textbf{59} (1984), no.~2, 335--395.

\bibitem{DP86}
Ben de~Pagter, \emph{Irreducible compact operators}, Mathematische Zeitschrift
  \textbf{192} (1986), no.~1, 149--153.

\bibitem{D10}
Klaus Deimling, \emph{Nonlinear functional analysis}, Courier Corporation,
  2010.

\bibitem{EGK20}
Katharina Eichinger, Manuel~V Gnann, and Christian Kuehn, \emph{Multiscale
  analysis for traveling-pulse solutions to the stochastic fitzhugh-nagumo
  equations}, arXiv preprint arXiv:2002.07234 (2020).

\bibitem{ENB00}
Klaus-Jochen Engel, Rainer Nagel, and Simon Brendle, \emph{One-parameter
  semigroups for linear evolution equations}, vol. 194, Springer, 2000.

\bibitem{E72i}
John~W Evans, \emph{Nerve axon equations: 1 linear approximations}, Indiana
  University Mathematics Journal \textbf{21} (1972), no.~9, 877--885.

\bibitem{F02}
Flavio~H Fenton, Elizabeth~M Cherry, Harold~M Hastings, and Steven~J Evans,
  \emph{Multiple mechanisms of spiral wave breakup in a model of cardiac
  electrical activity}, Chaos: An Interdisciplinary Journal of Nonlinear
  Science \textbf{12} (2002), no.~3, 852--892.

\bibitem{F58}
Ronald~Aylmer Fisher, \emph{The genetical theory of natural selection}, Ripol
  Classic, 1958.

\bibitem{FH61}
Richard FitzHugh, \emph{Impulses and physiological states in theoretical models
  of nerve membrane}, Biophysical journal \textbf{1} (1961), no.~6, 445--466.

\bibitem{FW98}
Mark~Iosifovich Freidlin and Alexander~D Wentzell, \emph{Random perturbations},
  Random perturbations of dynamical systems, Springer, 1998, pp.~15--43.

\bibitem{GMV20}
Alexandre Gaudilli{\`e}re, Paolo Milanesi, and Maria~Eul{\'a}lia Vares,
  \emph{Asymptotic exponential law for the transition time to equilibrium of
  the metastable kinetic ising model with vanishing magnetic field}, Journal of
  Statistical Physics \textbf{179} (2020), 263--308.

\bibitem{GV21}
Enrico Gavagnin, Sean~T Vittadello, Gency Gunasingh, Nikolas~K Haass, Matthew~J
  Simpson, Tim Rogers, and Christian~A Yates, \emph{Synchronized oscillations
  in growing cell populations are explained by demographic noise}, Biophysical
  journal \textbf{120} (2021), no.~8, 1314--1322.

\bibitem{GPS18}
Giambattista Giacomin, Christophe Poquet, and Assaf Shapira, \emph{Small noise
  and long time phase diffusion in stochastic limit cycle oscillators}, Journal
  of Differential Equations \textbf{264} (2018), no.~2, 1019--1049.

\bibitem{GM05}
Beniamin Goldys and Bohdan Maslowski, \emph{Exponential ergodicity for
  stochastic reaction-diffusion equations}, Stochastic Partial Differential
  Equations and Applications -- VII (2005), 115.

\bibitem{H13}
Martin Hairer, \emph{Solving the kpz equation}, Annals of mathematics (2013),
  559--664.

\bibitem{H14}
\bysame, \emph{A theory of regularity structures}, Inventiones mathematicae
  \textbf{198} (2014), no.~2, 269--504.

\bibitem{H20}
Christian Hendrik~Severian Hamster, \emph{Noisy patterns: Bridging the gap
  between stochastics and dynamics}, Ph.D. thesis, Leiden University, 2020.

\bibitem{HH20}
Christian Hendrik~Severian Hamster and Hermen~Jan Hupkes, \emph{Travelling
  waves for reaction--diffusion equations forced by translation invariant
  noise}, Physica D: Nonlinear Phenomena \textbf{401} (2020), 132233.

\bibitem{HK18}
Alexandru Hening and Martin Kolb, \emph{Quasistationary distributions for
  one-dimensional diffusions with singular boundary points}, Stochastic
  Processes and their Applications \textbf{129} (2019), no.~5, 1659--1696.

\bibitem{H22}
Alexandru Hening, Weiwei Qi, Zhongwei Shen, and Yingfei Yi, \emph{Population
  dynamics under demographic and environmental stochasticity}, arXiv preprint
  arXiv:2207.08883 (2022).

\bibitem{IM16}
James Inglis and James MacLaurin, \emph{A general framework for stochastic
  traveling waves and patterns, with application to neural field equations},
  SIAM Journal on Applied Dynamical Systems \textbf{15} (2016), no.~1,
  195--234.

\bibitem{JQSY22}
Min Ji, Weiwei Qi, Zhongwei Shen, and Yingfei Yi, \emph{Transient dynamics of
  absorbed singular diffusions}, Journal of Dynamics and Differential Equations
  \textbf{34} (2022), no.~4, 3089--3129.

\bibitem{KS08}
Panki Kim and Renming Song, \emph{Intrinsic ultracontractivity of non-symmetric
  diffusion semigroups in bounded domains}, Tohoku Mathematical Journal, Second
  Series \textbf{60} (2008), no.~4, 527--547.

\bibitem{KS14}
Jennifer Kr{\"u}ger and Wilhelm Stannat, \emph{Front propagation in stochastic
  neural fields: a rigorous mathematical framework}, SIAM Journal on Applied
  Dynamical Systems \textbf{13} (2014), no.~3, 1293--1310.

\bibitem{KS17}
\bysame, \emph{A multiscale-analysis of stochastic bistable reaction--diffusion
  equations}, Nonlinear Analysis \textbf{162} (2017), 197--223.

\bibitem{KSS97}
Martin Krupa, Bj{\"o}rn Sandstede, and Peter Szmolyan, \emph{Fast and slow
  waves in the fitzhugh--nagumo equation}, Journal of Differential Equations
  \textbf{133} (1997), no.~1, 49--97.

\bibitem{KMZ22}
Christian Kuehn, James MacLaurin, and Giulio Zucal, \emph{Stochastic rotating
  waves}, arXiv preprint arXiv:2111.07096 (2021).

\bibitem{L16}
Eva Lang, \emph{A multiscale analysis of traveling waves in stochastic neural
  fields}, SIAM Journal on Applied Dynamical Systems \textbf{15} (2016), no.~3,
  1581--1614.

\bibitem{LS16}
Eva Lang and Wilhelm Stannat, \emph{L2-stability of traveling wave solutions to
  nonlocal evolution equations}, Journal of Differential Equations \textbf{261}
  (2016), no.~8, 4275--4297.

\bibitem{LRR21}
Tony Leli{\`e}vre, Mouad Ramil, and Julien Reygner, \emph{Quasi-stationary
  distribution for the langevin process in cylindrical domains, part i:
  existence, uniqueness and long-time convergence}, Stochastic Processes and
  their Applications \textbf{144} (2022), 173--201.

\bibitem{LRSS21}
Rongli Liu, Yan-Xia Ren, Renming Song, and Zhenyao Sun, \emph{Quasi-stationary
  distributions for subcritical superprocesses}, Stochastic Processes and their
  Applications \textbf{132} (2021), 108--134.

\bibitem{L21}
Stefan Luther, Flavio~H Fenton, Bruce~G Kornreich, Amgad Squires, Philip
  Bittihn, Daniel Hornung, Markus Zabel, James Flanders, Andrea Gladuli, Luis
  Campoy, et~al., \emph{Low-energy control of electrical turbulence in the
  heart}, Nature \textbf{475} (2011), no.~7355, 235--239.

\bibitem{M20}
James MacLaurin, \emph{Metastability of waves and patterns subject to
  spatially-extended noise}, arXiv preprint arXiv:2006.12627 (2020).

\bibitem{M93}
Bohdan Maslowski, \emph{On probability distributions of solutions of semilinear
  stochastic evolution equations}, Stochastics: An International Journal of
  Probability and Stochastic Processes \textbf{45} (1993), no.~1-2, 17--44.

\bibitem{MS99}
Bohdan Maslowski and Jan Seidler, \emph{On sequentially weakly feller solutions
  to spde’s}, Atti della Accademia Nazionale dei Lincei. Classe di Scienze
  Fisiche, Matematiche e Naturali. Rendiconti Lincei. Matematica e Applicazioni
  \textbf{10} (1999), no.~2, 69--78.

\bibitem{MR81}
Robert~N. Miller and John Rinzel, \emph{The dependence of impulse propagation
  speed on firing frequency, dispersion, for the hodgkin-huxley model},
  Biophysical Journal \textbf{34} (1981), no.~2, 227--259.

\bibitem{MMQ11}
Carl Mueller, Leonid Mytnik, and Jeremy Quastel, \emph{Effect of noise on front
  propagation in reaction-diffusion equations of kpp type}, Inventiones
  mathematicae \textbf{184} (2011), no.~2, 405--453.

\bibitem{MMR21}
Carl Mueller, Leonid Mytnik, and Lenya Ryzhik, \emph{The speed of a random
  front for stochastic reaction--diffusion equations with strong noise},
  Communications in Mathematical Physics \textbf{384} (2021), no.~2, 699--732.

\bibitem{P12}
Amnon Pazy, \emph{Semigroups of linear operators and applications to partial
  differential equations}, vol.~44, Springer Science \& Business Media, 2012.

\bibitem{PZ95}
Szymon Peszat and Jerzy Zabczyk, \emph{Strong feller property and
  irreducibility for diffusions on hilbert spaces}, The Annals of Probability
  (1995), 157--172.

\bibitem{P85}
Ross~G Pinsky, \emph{On the convergence of diffusion processes conditioned to
  remain in a bounded region for large time to limiting positive recurrent
  diffusion processes}, The Annals of Probability (1985), 363--378.

\bibitem{P90}
\bysame, \emph{The lifetimes of conditioned diffusion processes}, Annales de
  l'IHP Probabilit{\'e}s et statistiques, vol.~26, 1990, pp.~87--99.

\bibitem{RMF08}
Tobias Reichenbach, Mauro Mobilia, and Erwin Frey, \emph{Self-organization of
  mobile populations in cyclic competition}, Journal of Theoretical Biology
  \textbf{254} (2008), no.~2, 368--383.

\bibitem{RW00}
Leonard~CG Rogers and David Williams, \emph{Diffusions, markov processes, and
  martingales: Volume 1, foundations}, vol.~1, Cambridge university press,
  2000.

\bibitem{R12}
Tim Rogers, Alan~J McKane, and Axel~G Rossberg, \emph{Demographic noise can
  lead to the spontaneous formation of species}, Europhysics Letters
  \textbf{97} (2012), no.~4, 40008.

\bibitem{RVY06}
Bernard Roynette, Pierre Vallois, and Marc Yor, \emph{Some penalisations of the
  wiener measure}, Japanese Journal of Mathematics \textbf{1} (2006), no.~1,
  263--290.

\bibitem{SBD19}
Michael Salins, Amarjit Budhiraja, and Paul Dupuis, \emph{Uniform large
  deviation principles for banach space valued stochastic evolution equations},
  Transactions of the American Mathematical Society \textbf{372} (2019),
  no.~12, 8363--8421.

\bibitem{SSW97}
Bj{\"o}rn Sandstede, Arnd Scheel, and Claudia Wulff, \emph{Dynamics of spiral
  waves on unbounded domains using center-manifold reductions}, journal of
  differential equations \textbf{141} (1997), no.~1, 122--149.

\bibitem{S47}
Henry Scheff{\'e}, \emph{A useful convergence theorem for probability
  distributions}, The Annals of Mathematical Statistics \textbf{18} (1947),
  no.~3, 434--438.

\bibitem{SS08}
Jonathan~A Sherratt and Matthew~J Smith, \emph{Periodic travelling waves in
  cyclic populations: field studies and reaction--diffusion models}, Journal of
  the Royal Society Interface \textbf{5} (2008), no.~22, 483--505.

\bibitem{S13}
Wilhelm Stannat, \emph{Stability of travelling waves in stochastic nagumo
  equations}, arXiv preprint arXiv:1301.6378 (2013).

\bibitem{vZ00}
Harry van Zanten, \emph{A multivariate central limit theorem for continuous
  local martingales}, Statistics \& probability letters \textbf{50} (2000),
  no.~3, 229--235.

\bibitem{W68}
Shinzo Watanabe, \emph{A limit theorem of branching processes and continuous
  state branching processes}, Journal of Mathematics of Kyoto University
  \textbf{8} (1968), no.~1, 141--167.

\bibitem{W74}
Arthur~T Winfree, \emph{Patterns of phase compromise in biological cycles},
  Journal of Mathematical Biology \textbf{1} (1974), no.~1, 73--93.

\bibitem{X93}
Jack Xin, \emph{Existence of a class of symmetric rotating spiral waves on
  finite disc domains in excitable media}, Indiana University Mathematics
  Journal (1993), 1305--1337.

\bibitem{Y96}
K{\"o}saku Yosida, \emph{Functional analysis}, Springer Science \& Business
  Media, 2012.

\bibitem{Z12}
Eberhard Zeidler, \emph{Applied functional analysis: main principles and their
  applications}, vol. 109, Springer Science \& Business Media, 2012.

\bibitem{ZLS14}
Junfei Zhang, Shoumei Li, and Renming Song, \emph{Quasi-stationarity and
  quasi-ergodicity of general markov processes}, Science China Mathematics
  \textbf{57} (2014), no.~10, 2013--2024.

\end{thebibliography}

\subsubsection*{Acknowledgements}
Thanks to J{\"u}rgen Jost, for his continuing patience and support, to H.L.~Duc, for his helpful advice, and several anonymous reviewers, for their immeasurably helpful criticism.  
his work was funded by the IMPRS, 

\end{document}